\documentclass[a4paper,12pt]{amsart}
\usepackage{amscd, amssymb, pgf, tikz, hyperref}

\newcommand{\morph}[2]{$#1$--$#2$ morph}
\newcommand{\LGmph}{\morph{\Lambda}{\Gamma}}
\newcommand{\subgrph}[2]{#1^{#2}}
\newcommand{\regsymbol}{\maltese}
\newcommand{\rfMm}[1]{\Mm^{\regsymbol}_{#1}}
\newcommand{\malteseref}{\eqref{eq:regmorph}}

\def\cs#1{{
\ensuremath{\mathbin{\stackrel{#1}{\leftharpoondown}}}}}

\newcommand{\id}{\operatorname{id}}
\newcommand{\dom}{\operatorname{dom}}
\newcommand{\cod}{\operatorname{cod}}
\newcommand{\Obj}{\operatorname{Obj}}
\newcommand{\Hom}{\operatorname{Hom}}
\newcommand{\Cov}{\operatorname{Cov}}

\newcommand{\End}{\operatorname{End}}

\newcommand{\field}[1]{\mathbb{#1}}
\newcommand{\CC}{\field{C}}
\newcommand{\NN}{\field{N}}
\newcommand{\TT}{\field{T}}
\newcommand{\ZZ}{\field{Z}}
\newcommand{\Cc}{\mathcal{C}}
\newcommand{\Hh}{\mathcal{H}}
\newcommand{\Kk}{\mathcal{K}}
\newcommand{\Ll}{\mathcal{L}}
\newcommand{\Mm}{\mathcal{M}}
\newcommand{\Oo}{\mathcal{O}}
\newcommand{\cocycle}{c}
\newcommand{\opX}{\mathop{X}}

\newcounter{mycount}
\newcounter{mySaveEnumi}

\newcommand{\tgrphlim}{{
    \renewcommand{\leftarrow}{\leftharpoondown}
    \varprojlim
}}

\def\clsp{\overline{\operatorname{span}}}

\newtheorem{thm}{Theorem}[section]
\newtheorem{cor}[thm]{Corollary}
\newtheorem{lem}[thm]{Lemma}
\newtheorem{prop}[thm]{Proposition}

\theoremstyle{definition}
\newtheorem{dfn}[thm]{Definition}

\theoremstyle{remark}
\newtheorem{rmk}[thm]{Remark}
\newtheorem{rmks}[thm]{Remarks}
\newtheorem{example}[thm]{Example}
\newtheorem{examples}[thm]{Examples}
\newtheorem*{examples*}{Examples}

\newtheorem{ntn}[thm]{Notation}

\numberwithin{equation}{section}

\title[$k$-morphs]{Generalised morphisms of $k$-graphs: $k$-morphs}

\author{Alex Kumjian}
\address{Alex Kumjian\\ Department of Mathematics (084)\\ University
of Nevada\\ Reno NV 89557-0084\\ USA} \email{alex@unr.edu}

\author{David Pask}
\address{David Pask, Aidan Sims\\ School of Mathematics and
Applied Statistics  \\
University of Wollongong\\
NSW  2522\\
AUSTRALIA} \email{dpask, asims@uow.edu.au}
\author{Aidan Sims}



\thanks{This research was supported by the Australian Research Council.}

\date{\today}

\begin{document}

\begin{abstract}
In a number of recent papers, $(k+l)$-graphs have been
constructed from $k$-graphs by inserting new edges in the last
$l$ dimensions. These constructions have been motivated by
$C^*$-algebraic considerations, so they have not been treated
systematically at the level of higher-rank graphs themselves.
Here we introduce $k$-morphs, which provide a systematic
unifying framework for these various constructions. We think of
$k$-morphs as the analogue, at the level of $k$-graphs, of
$C^*$-correspondences between $C^*$-algebras. To make this
analogy explicit, we introduce a category whose objects are
$k$-graphs and whose morphisms are isomorphism classes of
$k$-morphs. We show how to extend the assignment $\Lambda \mapsto
C^*(\Lambda)$ to a functor from this category to the category
whose objects are $C^*$-algebras and whose morphisms are
isomorphism classes of $C^*$-correspondences.
\end{abstract}

\maketitle

\section{Introduction}
Over the last ten years, graph $C^*$-algebras and their
analogues have been the subject of intense research interest
(see for example \cite{Dea:continuous graphs, DHS, DT, FLR, K4,
KPRR, MT, Spiel:functorial}, or see \cite{CBMSbk} for a good
overview). In particular, the higher-rank graphs and associated
$C^*$-algebras introduced in~\cite{KP} have recently been
widely studied \cite{Evans, FMY, KP3, PRRS}. Higher-rank graphs
generalise directed graphs, so there are many points of
similarity between the two theories, especially at the level of
fundamental existence and uniqueness results. However, as both
fields progress, the two sets of results are diverging more and
more rapidly.

One reason for this is the relatively involved combinatorial
structure of higher-rank graphs as opposed to ``ordinary"
one-dimensional graphs. It is fairly straightforward to modify
an ordinary graph by simply adding vertices or edges because
these are local operations. By contrast, adding vertices and
edges to a higher-rank graph is quite complicated because the
combinatorial peculiarities of higher-rank graphs mean that the
addition of an edge at some vertex typically necessitates
similar changes throughout a large portion of the higher-rank
graph. A good illustration of this is the contrast between the
straightforward process of ``adding tails" to a directed graph
\cite{BPRS} and the analogous but vastly more complicated
``removing sources" construction for higher-rank graphs
\cite{Far}.

It has become clear recently, however, that if higher-rank
graphs are not well-suited to constructions which involve
localised modifications, they are amenable to a somewhat
different style of construction which is not available in the
one-dimensional setting and which is proving very profitable
from a $C^*$-algebraic standpoint. Specifically, $k$-graphs
lend themselves to constructions whereby one increases the rank
of a graph or graphs by adding edges in new dimensions
\cite{FPS,KPS, PQS}. The resulting $(k+l)$-graph $C^*$-algebras
have been analysed as direct limits \cite{KPS, PQS} and as
crossed-products by group actions \cite{FPS}.

So far these constructions of $(k+l)$-graphs from $k$-graphs
have been ad hoc: in each of \cite{FPS,KPS, PQS}, given a
$k$-graph with natural symmetry or a pair of $k$-graphs with
structural similarities, the authors have constructed a
$(k+l)$-graph with bare hands. In each case, the $(k+l)$-graph
contains copies of the original $k$-graph or -graphs in the
first $k$-dimensions, and encodes the additional symmetry or
structural similarities in the remaining $l$ dimensions.

The purpose of this article is to replace these ad hoc methods
with a unifying construction which is functorial with respect
to the assignment of $C^*$-algebras to higher-rank graphs. More
specifically, in Section~\ref{sec:morphs} we axiomatise the
data required to insert a set $X$ of edges in a $(k+1)^{\rm
st}$ dimension between vertices in a $k$-graph $\Gamma$ and
those in a $k$-graph $\Lambda$ so as to obtain a $(k+1)$-graph.
We call a set $X$ endowed with such data a \LGmph, or a
$k$-morph from $\Gamma$ to $\Lambda$. Given a
\morph{\Lambda_0}{\Lambda_1}\ $X_1$ and a
\morph{\Lambda_1}{\Lambda_2}\ $X_2$, we define a fibred product
$X_1 *_{\Lambda_1^0} X_2$ which is a
\morph{\Lambda_0}{\Lambda_2}.
%% We show the isomorphism class of $X_1 *_{\Lambda_1^0} X_2$
%% depends only on the isomorphism classes of $X_1$ and $X_2$.
We show in Theorem~\ref{thm:category} that there is a category
$\Mm_k$ whose objects are $k$-graphs and whose morphisms are
isomorphism classes of $k$-morphs.

In Sections \ref{sec:linking}~and~\ref{sec:systems} we discuss
how $k$-morphs can be used as a model for constructions such as
those of \cite{FPS,KPS, PQS}. Given $k$-graphs $\Lambda$ and
$\Gamma$, and a \LGmph\ $X$, we define in
Section~\ref{sec:linking} what we call a \emph{linking graph}
for $X$. Roughly speaking, a linking graph for $X$ is a
$(k+1)$-graph $\Sigma$ containing disjoint copies of $\Lambda$
and $\Gamma$ connected in the $(k+1)^{\rm st}$ dimension by a
copy of $X$. We show in Proposition~\ref{prp:linking graph
exists} that a linking graph always exists, is unique up to
isomorphism, and is determined up to isomorphism by the
isomorphism class of $X$.

The constructions set out in \cite{FPS, KPS, PQS} typically
involve a system of linking graphs which are glued together in
some systematic way. For example, we can think of the covering
systems of \cite{KPS} as a system, organised by an underlying
Bratteli diagram, of linking graphs for $k$-morphs determined
by $k$-graph coverings. To capture this idea, we introduce in
Section~\ref{sec:systems} the notion of a $\Gamma$-system of
$k$-morphs, and the notion of a $\Gamma$-bundle for a
$\Gamma$-system. Given an $l$-graph $\Gamma$, a $\Gamma$ system
consists of a collection $\{\Lambda_v : v \in \Gamma^0\}$ of
$k$-graphs connected by $k$-morphs $\{X_\gamma : \gamma \in
\Gamma\}$ so that composition in $\Gamma$ corresponds in a
consistent way to the fibred product operation on the
associated $k$-morphs. A $\Gamma$-bundle for this system is
then a $(k+l)$-graph $\Sigma$ together with a map $f : \Sigma
\to \Gamma$ such that $f^{-1}(v) \cong \Lambda_v$ for each $v
\in \Gamma^0$ and such that $f^{-1}(\{r(\gamma), \gamma,
s(\gamma)\})$ is a linking graph for $X_\gamma$ for each
$\gamma \in \Gamma$. We call the map $f$ the bundle map for the
$\Gamma$-bundle $\Sigma$. We show in
Theorem~\ref{thm:skewgraphexists} that every $\Gamma$-system
admits a $\Gamma$-bundle, and that the $\Gamma$-bundle is
unique up to isomorphism and depends only on the isomorphism
class of the $\Gamma$-system. We indicate how to realise the
$k$-graphs constructed in \cite{FPS, KPS, PQS, PRRS} as
$\Gamma$-bundles in a natural way.

A $\Gamma$-system $X$ of $k$-morphs determines a functor from
$\Gamma$ into the category $\Mm_k$ via the assignments $v
\mapsto \Lambda_v$ and $\gamma \mapsto [X_\gamma]$ (where
$[X_\gamma]$ is the isomorphism class of $X_\gamma$). One might
initially hope that the isomorphism class of a $\Gamma$-bundle
for the system would be determined by this functor, so that we
could replace $\Gamma$-systems with functors. We show in
Proposition \ref{prp:1-graph systems} that when $\Gamma$ is a
$1$-graph each $\Gamma$-system is indeed determined up to
isomorphism by the functor $\gamma \mapsto [X_\gamma]$. However
this is the best we can hope for: Example~\ref{egs:not just
functors}(\ref{it:non-unique system}) shows that if $\Gamma$
has rank 2 or more there may be non-isomorphic $\Gamma$-systems
which determine the same functor from $\Gamma$ to $\Mm_k$; and
an example of Spielberg's, which we present as
Example~\ref{egs:not just functors}(\ref{it:Jacks}), shows that
there exists a 3-graph $\Gamma$ and a functor $\Gamma \to
\Mm_0$ which is not the functor determined by any
$\Gamma$-system of $0$-graphs. In particular, $\Gamma$-systems
cannot be replaced with functors.

An example of a $\Gamma$ system is the following. Let $X$ be a
\morph{\Lambda}{\Lambda}\ (we refer to these as $\Lambda$
endomorphs). Then $X$ gives rise to a $T_1$-system of
$k$-morphs, where $T_1$ is the $1$-graph with a single vertex
and a single edge. We call the $T_1$-bundle for such a system
the endomorph skew graph of $\Lambda$ by $X$, and denote it
$\Lambda \times_X \NN$. When $X$ arises from an automorphism of
$\Lambda$, we recover the crossed product graph of \cite{FPS}.

In Section~\ref{sec:functoriality} we discuss how our
constructions behave at the level of $C^*$-algebras. The
category $\Mm_k$ is reminiscent of the category (which we shall
denote $\Cc$) of \cite{EKQR, Landsman:bicategories, Sch}, whose
objects are $C^*$-algebras and whose morphisms are isomorphism
classes of $C^*$-correspondences (also known as Hilbert
bimodules). To simplify arguments, we restrict our attention to
a subcategory $\rfMm{k}$ of $\Mm_k$. We construct a
$C^*$-correspondence $\Hh(X)$ for each $k$-morph $X$, in such a
way that the isomorphism class of $\Hh(X)$ depends only on that
of $X$. We show in Theorem~\ref{thm:functoriality} that the
assignments $\Lambda \mapsto C^*(\Lambda)$ and $[X] \mapsto
[\Hh(X)]$ determine a contravariant functor between $\rfMm{k}$
and $\Cc$. In the special case where $X$ is a
$\Lambda$-endomorph, so that $[X] \in
\End_{\rfMm{k}}(\Lambda)$, Theorem~\ref{thm:CP algebra} shows
that the $C^*$-algebra of the endomorph
skew graph % $T_1$-bundle for $X$
is canonically isomorphic to the Cuntz-Pimsner algebra
$\Oo_{\Hh(X)}$.

\vskip0.5em

\textbf{Acknowledgements.} The authors wish to thank Jack
Spielberg for providing the as yet unpublished example \cite{JS
example} which we have reproduced in Examples~\ref{egs:not just
functors}\eqref{it:Jacks}. Much of this work was completed
during the first author's recent trip to Australia. He wishes
to thank his colleagues for their hospitality. We would also
like to acknowledge the support and hospitality of the Fields
Institute during the final stages of preparation of this
manuscript.

\section{Preliminaries}\label{sec:prelims}

\subsection{Higher-rank graphs}\label{sec:hrgs}
In this paper, unlike previous treatments of $k$-graphs
\cite{Evans, FMY, KP, RSY1}, we allow $0$-graphs. To make sense
of this, we take the convention that $\NN^0$ is the trivial
semigroup $\{0\}$. We will also insist that all $k$-graphs are
nonempty.

Modulo the minor differences mentioned above, we will adopt the
conventions of \cite{KP, PQR} for $k$-graphs. Given a
nonnegative integer $k$, a \emph{$k$-graph} is a nonempty
countable small category $\Lambda$ equipped with a functor $d
:\Lambda \to \NN^k$ satisfying the \emph{factorisation
property}: for all $\lambda \in \Lambda$ and $m,n \in \NN^k$
such that $d( \lambda )=m+n$ there exist unique $\mu ,\nu \in
\Lambda$ such that $d(\mu)=m$, $d(\nu)=n$, and $\lambda=\mu
\nu$. When $d(\lambda )=n$ we say $\lambda$ has \emph{degree}
$n$. By abuse of notation, we will use $d$ to denote the degree
functor in every $k$-graph in this paper; the domain of $d$ is
always clear from context.

For $k \ge 1$, the standard generators of $\NN^k$ are denoted $e_1,
\dots, e_k$, and for $n \in \NN^k$ and $1 \le i \le k$ we write $n_i$
for the $i^{\rm th}$ coordinate of $n$.

For $n \in \NN^k$, we write $\Lambda^n$ for $d^{-1}(n)$. In
particular, $\Lambda^0$ is the vertex set. The \emph{vertices}
of $\Lambda$ are the elements of $\Lambda^0$. The factorisation
property implies that $o \mapsto \id_o$ is a bijection from the
objects of $\Lambda$ to $\Lambda^0$. We will frequently use
this bijection to silently identify $\Obj(\Lambda)$ with
$\Lambda^0$. The domain and codomain maps in the category
$\Lambda$ therefore become maps $s,r : \Lambda \to \Lambda^0$.
More precisely, for $\alpha \in\Lambda$, the \emph{source}
$s(\alpha)$ is the identity morphism associated with the object
$\dom(\alpha)$ and similarly, $r(\alpha) = \id_{\cod(\alpha)}$.

Note that a $0$-graph is then a countable category whose only
morphisms are the identity morphisms; we think of them as a
collection of isolated vertices.

For $u,v\in\Lambda^0$ and $E \subset \Lambda$, we write $u E$
for $E \cap r^{-1}(u)$ and $E v$ for $E \cap s^{-1}(v)$. For $n
\in \NN^k$, we denote by $\Lambda^{\le n}$ the set
\[
\Lambda^{\le n} = \{\lambda \in \Lambda : d(\lambda) \le n , s(\lambda)
\Lambda^{e_i}  = \emptyset \text{ whenever } d(\lambda) + e_i \le n\} .
\]

We say that $\Lambda$ is \emph{row-finite} if $v\Lambda^n$ is finite
for all $v \in \Lambda^0$ and $n \in \NN^k$. We say that $\Lambda$ is
\emph{locally convex} if whenever $1 \le i < j \le k$, $e \in
\Lambda^{e_i}$, $f \in \Lambda^{e_j}$ and $r(e) = r(f)$, we can
extend both $e$ and $f$ to paths $ee'$ and $ff'$ in $\Lambda^{e_i +
e_j}$.

\subsection{Maps between higher-rank graphs}\label{sec:maps}
A \emph{$k$-graph morphism} is a degree-preserving functor.
More generally, if $\omega : \NN^k \to \NN^l$ is a
homomorphism, $\Lambda$ is a $k$-graph and $\Gamma$ is an
$l$-graph, we say that a functor $f : \Lambda \to \Gamma$ is an
\emph{$\omega$-quasimorphism} if $d_\Gamma(f(\lambda)) =
\omega(d_\Lambda(\lambda))$ for all $\lambda \in \Lambda$. A
$k$-graph morphism is then an $\id_k$-quasimorphism.

Let $\omega : \NN^k \to \NN^l$ be a homomorphism, and let
$\Gamma$ be an $l$-graph. The pullback $\omega^*\Gamma$ is the
$k$-graph $\omega^*\Gamma = \{(\gamma,n) \in \Gamma \times
\NN^k : \omega(n) = d(\gamma)\}$ with degree map
$d_{\omega^*\Gamma}(\gamma,n) = n$ \cite[Definition~1.9]{KP}.
In the case where $\omega$ is injective, it will also sometimes
be convenient to regard the subcategory $\Gamma^{\omega} :=
\bigcup_{n \in \NN^k} \Gamma^{\omega(n)}$ of $\Gamma$ as a
$k$-graph as follows. We define the degree functor $d_\omega$
on $\Gamma^\omega$ by $d_{\omega}(\gamma) = n$ when $\omega(n)
= d_\Gamma(\gamma)$. Of course $\Gamma^\omega$ and
$\omega^*\Gamma$ are isomorphic, but the former is a subset of
$\Gamma$ whereas the latter is formally disjoint from $\Gamma$.

As in \cite{PQR}, a \emph{covering} of a $k$-graph $\Lambda$ by
a $k$-graph $\Gamma$ is a surjective $k$-graph morphism $p
:\Gamma \to \Lambda$ such that for all $v\in \Gamma^0$, $p$
restricts to bijections between $v\Gamma$ and $p(v)\Lambda$ and
between $\Gamma v$ and $\Lambda p(v)$. The covering $p : \Gamma
\to \Lambda$ is \emph{finite} if $p^{-1} (v)$ is finite for all
$v \in \Lambda^0$. Every covering $p : \Gamma \to \Lambda$ has
the unique path lifting property: for every $\lambda \in
\Lambda$ and $v \in \Gamma^0$ with $p(v)=s(\lambda)$ there is a
unique $\gamma \in \Gamma$ such that $p(\gamma) = \lambda$ and
$s(\gamma) = v$; and similarly at $r(\lambda)$.

\subsection{\texorpdfstring{$C^*$}{C*}-algebras associated to higher-rank
graphs}\label{sec:intro C*}

Given a row-finite, locally convex $k$-graph $(\Lambda, d)$, a
Cuntz-Krieger $\Lambda$-family is a collection $\{t_\lambda :
\lambda \in \Lambda \}$ of partial isometries satisfying the
Cuntz-Krieger relations:
\begin{itemize}
\item $\{t_v : v \in \Lambda^0\}$ is a collection of
    mutually orthogonal projections;
\item $t_\lambda t_\mu = t_{\lambda \mu}$ whenever
    $s(\lambda) = r(\mu)$;
\item $t^*_\lambda t_\lambda = t_{s(\lambda)}$ for all
    $\lambda \in \Lambda$; and
\item $t_v = \sum_{\lambda \in v \Lambda^{\le n}} t_\lambda
    t^*_\lambda$ for all $v \in \Lambda^0$ and $n \in
    \NN^k$.
\end{itemize}
The $k$-graph $C^*$-algebra $C^*(\Lambda)$ is the universal
$C^*$-algebra generated by a Cuntz-Krieger $\Lambda$-family
$\{s_\lambda : \lambda \in \Lambda \}$. That is, for every
Cuntz-Krieger $\Lambda$-family $\{t_\lambda : \lambda \in
\Lambda \}$ there is a homomorphism $\pi_t$ of $C^*( \Lambda )$
satisfying $\pi_t(s_\lambda) = t_\lambda$ for all $\lambda \in
\Lambda$.

A $k$-graph with no sources is automatically locally convex
with $\Lambda^{\le n} = \Lambda^n$ for all $n \in \NN^k$. Hence
the definition of $C^*(\Lambda)$ above reduces in this case to
\cite[Definition 1.5]{KP}.

By \cite[Theorem 3.15]{RSY1}, the generating partial isometries
$\{s_\lambda : \lambda \in \Lambda\} \subset C^*(\Lambda)$ are
all nonzero.

If $\Lambda$ is a $0$-graph, then it trivially has no sources,
and the last three Cuntz-Krieger relations follow from the
first one. So $C^*(\Lambda)$ is the universal $C^*$-algebra
generated by mutually orthogonal projections $\{s_v : v \in
\Lambda^0\}$; that is $C^*(\Lambda) \cong c_0(\Lambda^0)$.

Let $\Lambda$ be a $k$-graph. There is a strongly continuous
action $\gamma$ of $\TT^k$ on $C^*(\Lambda)$, called the
\emph{gauge-action}, such that $\gamma_z(s_\lambda) =
z^{d(\lambda)} s_\lambda$ for all $z \in \TT^k$ and $\lambda
\in \Lambda$.

\subsection{\texorpdfstring{$C^*$}{C*}-correspondences}\label{sec:correspondences}

We define Hilbert modules following \cite{Lan} and
\cite[\S II.7]{Black}.
Let $B$ be a $C^*$-algebra and let $\Hh$ be a right $B$-module.
Then a {\em $B$-valued inner product} on $\Hh$ is a function
$\langle \cdot, \cdot\rangle_B : \Hh \times  \Hh \to B$ satisfying the
following conditions for all $\xi, \eta, \zeta \in \Hh$, $b \in B$ and
$\alpha, \beta \in \CC$:
\begin{itemize}
\item
$\langle \xi, \alpha\eta + \beta\zeta\rangle_B =
\alpha\langle \xi, \eta\rangle_B + \beta\langle \xi, \zeta\rangle_B$,
\item
$\langle \xi, \eta b\rangle_B = \langle \xi, \eta\rangle_Bb$,
\item $\langle \xi, \eta\rangle_B = \langle \eta, \xi\rangle_B^*$,
\item  $\langle \xi, \xi\rangle_B \ge 0$  and $\langle \xi,
    \xi\rangle_B = 0$ if and only if $\xi = 0$.
\end{itemize}
If $\Hh$ is complete with respect to the norm given by $\| \xi
\|^2 = \langle \xi, \xi\rangle_B$ then $\Hh$ is said to be a
(right-) {\em Hilbert $B$-module}.  If the range of the inner
product is not contained in any proper ideal in $B$, $\Hh$ is
said to be {\em full}. Note that $B$ may be endowed with the
structure of a full Hilbert $B$-module by taking $\langle \xi,
\eta\rangle_B = \xi^*\eta$ for all $\xi, \eta \in B$. A map $T
: \Hh \to \Hh$ is an {\em adjointable operator} if there is a
map $T^* : \Hh \to \Hh$ such that $\langle T\xi, \eta\rangle_B
= \langle \xi, T^*\eta\rangle_B$ for all $\xi, \eta \in \Hh$.
Such an operator is necessarily linear and bounded and the
collection $\Ll(\Hh)$ of all adjointable operators on $\Hh$ is
a $C^*$-algebra. Each pair $\xi, \eta \in \Hh$ determines a
\emph{rank-one operator} $\theta_{\xi,\eta}$ (with adjoint
$\theta_{\eta, \xi}$) given by $\theta_{\xi, \eta}\zeta =
\xi\langle \eta, \zeta\rangle_B$ for $\zeta \in \Hh$. The
$C^*$-subalgebra $\Kk(\Hh)$ of $\Ll(\Hh)$ generated by the
$\theta_{\xi, \eta}$ is called the algebra of compact operators
on $\Hh$.  Note that $\Ll(\Hh)$ may be identified with the
multiplier algebra of $\Kk(\Hh)$.

Let $A$ and $B$ be $C^*$-algebras; then a $C^*$-correspondence
from $A$ to $B$ or more briefly an $A$--$B$
$C^*$-correspondence is a Hilbert $B$-module $\Hh$ together
with a $*$-homomorphism $\varphi : A \to \Ll(\Hh)$. Given a
homomorphism $\varphi: A \to B$, we may endow $B$ with the
structure of $A$--$B$ $C^*$-correspondence in a canonical way.
So it is natural to think of an $A$--$B$ $C^*$-correspondence
as a generalised homomorphism from $A$ to $B$. A
$C^*$-correspondence $\Hh$ is said to be {\em nondegenerate} if
span\,$\{ \varphi(a)\xi : a \in A, \xi \in \Hh \}$ is dense in
$\Hh$ (some authors have also called such $C^*$-correspondences
\emph{essential}). We often suppress $\varphi$ by writing $a
\cdot \xi$ for $\varphi(a)\xi$.

As discussed in \cite{Black,EKQR, Landsman:bicategories, Sch},
there is a category $\Cc$ with $\Obj(\Cc)$ the class of
$C^*$-algebras, and $\Hom_{\Cc}(A, B)$ the set of isomorphism
classes of $A$--$B$ $C^*$-correspondences with identity
morphisms $[A]$. Composition
\[
\Hom_{\Cc}(B, C) \times \Hom_{\Cc}(A, B) \to
\Hom_{\Cc}(A, C)
\]
is defined by $([\Hh_1],[\Hh_2]) \mapsto [\Hh_2 \otimes_B
\Hh_1]$ where $\Hh_2 \otimes_B \Hh_1$ denotes the tensor
product of $C^*$-correspondences. The $C^*$-correspondence
$\Hh_1 \otimes \Hh_2$ is called the \emph{internal tensor
product} of $\Hh_1$ and $\Hh_2$ by Blackadar and the
\emph{interior tensor product} by Lance (see
\cite[II.7.4.1]{Black} and \cite[Prop.~4.5]{Lan} and the
following discussion).

\subsection{Representations of \texorpdfstring{$C^*$}{C*}-correspondences}

Let $\Hh$ be an $A$--$A$ $C^*$-correspondence. Recall from
\cite{Pim} that a representation of $\Hh$ in a $C^*$-algebra
$B$ is a pair $(t, \pi)$ where $\pi : A \to B$ is a
homomorphism, $t : \Hh \to B$ is linear, and such that for all
$a \in A$ and $\xi, \eta \in \Hh$, we have $t(a\cdot \xi) =
\pi(a)t(\xi)$, $t(\xi \cdot a) = t(\xi)\pi(a)$, and
$\pi(\langle \xi, \eta \rangle_A) = t(\xi)^* t(\eta)$.

Given a $C^*$-correspondence $\Hh$ over $A$ and a
representation $(t, \pi)$ of $\Hh$ on $B$, there is a
homomorphism $t^{(1)} : \Kk(\Hh) \to B$ satisfying
$t^{(1)}(\theta_{\xi,\eta}) = t(\xi) t(\eta)^*$ for all
$\xi,\eta \in \Hh$ (Pimsner denotes this homomorphism
$\pi^{(1)}$ in \cite{Pim}). In the cases of interest later in
this paper, the left action of $A$ on $\Hh$ is by elements of
$\Kk(\Hh)$ (that is $\varphi : A \to \Ll(\Hh)$ in fact takes
values in $\Kk(\Hh)$), so $t^{(1)} \circ \varphi$ is a
homomorphism from $A$ to $B$. In this case, the pair $(t, \pi)$
is said to be \emph{Cuntz-Pimsner covariant} if $t^{(1)} \circ
\varphi = \pi$.

Given a $C^*$-correspondence $\Hh$ over $A$, there is a
representation $(j_\Hh, j_A)$ in a $C^*$-algebra $\Oo_\Hh$
which is universal in the sense that given another
representation $(t, \pi)$ of $\Hh$ in $B$ there is a
homomorphism $t \times \pi : \Oo_\Hh \to B$ satisfying $(t
\times \pi) \circ j_\Hh = t$ and $(t \times \pi) \circ j_A =
\pi$.

\section{\texorpdfstring{$k$}{k}-morphs}\label{sec:morphs}

In this section, we define $k$-morphs, provide some motivating
examples, and show how isomorphism classes of $k$-morphs can be
regarded as the morphisms of a category whose objects are
$k$-graphs. Conceptually, a $k$-morph may be thought of as a
bridge between two $k$-graphs $\Lambda$ and $\Gamma$; it
consists of a set $X$ and some structure maps which are
precisely what is needed to build a $(k+1)$-graph that contains
disjoint copies of $\Lambda$ and $\Gamma$ and in which elements
of $X$ become edges of degree $e_{k+1}$ from vertices in the
copy of $\Gamma$ to vertices in the copy of $\Lambda$. We now
give the formal definition.

\begin{dfn}\label{dfn:k-morphs}
Let $\Lambda$ and $\Gamma$ be $k$-graphs, let $X$ be a
countable set, and fix functions $r : X \to \Lambda^0$ and $s :
X \to \Gamma^0$. We will write $X *_{\Gamma^0} \Gamma$ for the
fibred product $\{(x,\gamma) : x \in X, \gamma \in \Gamma, s(x)
= r(\gamma)\}$. Likewise, we will write $\Lambda *_{\Lambda^0}
X$ for the fibred product $\{(\lambda,x) : \lambda \in \Lambda,
x \in X, s(\lambda) = r(x)\}$. Fix a bijection $\phi : X
*_{\Gamma^0} \Gamma \to \Lambda *_{\Lambda^0} X$, and suppose
that whenever $\phi(x_1, \gamma_1) = (\lambda_1, x_2)$, we have
\begin{enumerate}
    \item\label{it:morph d} $d(\gamma_1) = d(\lambda_1)$;
    \item\label{it:morph s} $s(\gamma_1) = s(x_2)$; and
    \item\label{it:morph r} $r(\lambda_1) =
        r(x_1)$.\setcounter{mySaveEnumi}{\value{enumi}}
\end{enumerate}
Suppose further that whenever $\phi(x_1, \gamma_1) =
(\lambda_1, x_2)$ and $\phi(x_2, \gamma_2) = (\lambda_2,
x_3)${\renewcommand{\thefootnote}{\dag}\footnote{Note that the
conditions $\phi(x_1, \gamma_1) = (\lambda_1, x_2)$ and
$\phi(x_2, \gamma_2) = (\lambda_2, x_3)$ together with
(2)~and~(3) imply that $s(\gamma_1) = s(x_2) = r(\gamma_2)$ and
$s(\lambda_1) = r(x_2) = r(\lambda_2)$; it then follows that
$(\lambda_1\lambda_2, x_3) \in \Lambda *_{\Lambda^0} X$ and
$(x_1, \gamma_1\gamma_2) \in X *_{\Gamma^0} \Gamma$.}}, we have
\begin{enumerate}\setcounter{enumi}{\value{mySaveEnumi}}
    \item\label{it:phi assoc} $\phi(x_1, \gamma_1\gamma_2)
        = (\lambda_1\lambda_2, x_3)$.
\end{enumerate}
Then we call $X$  a \emph{\LGmph}, or simply a $k$-morph. If
$\Lambda = \Gamma$, then we call $X$ a \emph{$\Lambda$
endomorph}.
\end{dfn}

\begin{rmk}
Technically a \LGmph\ is a quadruple $(X, r, s, \phi)$, but by
the usual abuse of notation, we will say ``$X$ is a \LGmph''
without reference to the additional structure.
\end{rmk}

\begin{examples} \label{ex:first}
We now present a series of examples of $k$-morphs. In each case
we shall describe the set $X$ and the structure maps; it is
straightforward to check in each case that the resulting data
define a $k$-morph.

\begin{enumerate}\renewcommand{\theenumi}{\roman{enumi}}
\item \label{it:X(alpha)}Let $\Lambda$ and $\Gamma$ be
    $k$-graphs, and let $\alpha : \Gamma \to \Lambda$ be an
    isomorphism. Let $X(\alpha) = \Gamma^0$ and define
    structure maps $r_\alpha = \alpha$, $s_\alpha =
    \id_{\Gamma^0}$ and $\phi(r(\gamma), \gamma) =
    (\alpha(\gamma), s(\gamma))$. Then $X(\alpha)$ is a
    \LGmph. If $\Lambda = \Gamma$ so that $\alpha$ is an
    automorphism, then $X(\alpha)$ is a $\Lambda$
    endomorph. In the special case where $\alpha$ is the
    identity isomorphism $\id_\Lambda$, we refer to
    $X(\id_\Lambda)$ as the \emph{identity endomorph} on
    $\Lambda$. When it is useful to highlight its
    dependence on $\Lambda$ we will denote it as
    $I_\Lambda$.
\item\label{it:pX} Let $p : \Gamma \to \Lambda$ be a
    covering map. Let $\sideset{_p}{}{\opX} = \Gamma^0$,
    and define structure maps by $r = p$, $s =
    \id_{\Gamma^0}$ and $\phi(r(\gamma), \gamma) =
    (p(\gamma), s(\gamma))$. Then $\sideset{_p}{}{\opX}$ is
    a \LGmph. Such a $k$-morph is called a \emph{covering
    $k$-morph}. Note that if $p = \alpha$ is an
    isomorphism, then $\sideset{_p}{}{\opX}$ is equal to
    the $k$-morph $X(\alpha)$ of the preceding example.
\item\label{it:Xq} We may reverse the ``direction'' of the
    elements of $X$ in the preceding example to get a
    \morph{\Gamma}{\Lambda}. Let $p : \Gamma \to \Lambda$
    be a covering of $k$-graphs. Let $\opX_p := \Gamma^0$,
    and define $r = \id_{\Gamma^0}$, $s = p$ and
    $\phi(r(\gamma), p(\gamma)) = (\gamma, s(\gamma))$
    (where we are using the unique path lifting property to
    recover $\gamma$ from $p(\gamma)$ and $r(\gamma)$).
    Then $\opX_p$ is a \morph{\Gamma}{\Lambda}.
\item \label{it:pXq} Let $\Lambda_1, \Lambda_2, \Gamma$ be
    $k$-graphs and $p : \Gamma \to \Lambda_1$, $q : \Gamma
    \to \Lambda_2$ be coverings. Let
    $\sideset{_p}{_q}{\opX} = \Gamma^0$ and define
    structure maps by $r = p$, $s = q$, and
    $\phi(r(\gamma), q(\gamma)) = (p(\gamma), s(\gamma))$.
    Then $\sideset{_p}{_q}{\opX}$ is a
    \morph{\Lambda_1}{\Lambda_2}. This generalises the
    preceding two examples: if $\Lambda_2 = \Gamma$ and $q
    = \id_\Gamma$, then $\sideset{_p}{_{\id_\Gamma}}{\opX}
    = \sideset{_p}{}{\opX}$, and similarly if $\Lambda_1 =
    \Gamma$, then $\sideset{_{\id_\Lambda}}{_q}{\opX} =
    X_q$.
\item Number~(\ref{it:pXq}) (hence also numbers
    (\ref{it:pX})~and~(\ref{it:Xq})) above can be enriched
    with multiple ``edges" as in \cite{KPS}. Let $p :
    \Gamma \to \Lambda_1$ and $q : \Gamma \to \Lambda_2$ be
    covering maps. Write $S_m$ for the group of
    permutations of $\{1, \dots, m\}$; let $\cocycle$ be a
    cocycle from $\Gamma$ to $S_m$ (that is
    $\cocycle(\alpha)\cocycle(\beta) =
    \cocycle(\alpha\beta)$ whenever $\alpha$ and $\beta$
    are composable in $\Gamma$). Set
    $\sideset{_p^\cocycle}{_q}{\opX} = \Gamma^0 \times \{1,
    \dots, m\}$, and define structure maps by $r(v,i) =
    p(v)$, $s(v,i) = q(v)$ and $\phi((r(\gamma), i),
    q(\gamma)) = (p(\gamma), (s(\gamma),
    \cocycle(\gamma)^{-1}i))$. Then
    $\sideset{_p^\cocycle}{_q}{\opX}$ is a
    \morph{\Lambda_1}{\Lambda_2}.
\item\label{it:Sigmaiota} Let $(\Sigma, d)$ be a
    $(k+1)$-graph and $\iota : \NN^{k} \to \NN^{k+1}$ be
    the homomorphism $n \to (n, 0)$. Recall from
    Section~\ref{sec:maps} that we can regard $\Sigma^\iota
    := \{\lambda \in \Sigma : d(\lambda) \in
    \iota(\NN^k)\}$ as a $k$-graph. Let $X =
    \Sigma^{e_{k+1}}$, and define $r_X, s_X : X \to
    \Sigma^0$ to be the range and source maps $r,s$
    inherited from $\Sigma$. The bijection $\phi$ is
    obtained from the factorisation property in $\Sigma$:
    $\phi(x,\lambda) = (\lambda',x')$ where $x' \in X$ and
    $\lambda' \in \Sigma^\iota$ are the unique elements
    satisfying $x\lambda = \lambda'x'$ in $\Sigma$. Then
    $X$ is a $\Sigma^\iota$ endomorph.
\item Let $\Lambda$ and $\Gamma$ be $k$-graphs, and let
    $X_1$ and $X_2$ be {\LGmph}s. Then $X := X_1 \sqcup
    X_2$ is a \LGmph\ with the inherited structure maps.
\end{enumerate}
\end{examples}

We next define a kind of fibred product of $k$-morphs. This
fibred product, like tensor products, is not quite associative
on $k$-morphs. However Proposition~\ref{prp:composition} and
Lemma~\ref{lem:assoc} show that it does determine an
associative binary operation on isomorphism classes of
$k$-morphs. Of course, we must first say exactly what we mean
by an isomorphism of $k$-morphs.

\begin{dfn}
Fix $k$-graphs $\Lambda$ and $\Gamma$. Let $X$ and $Y$ be
{\LGmph}s. We say that $X$ and $Y$ are \emph{isomorphic} if
there is a bijection $\theta : X \to Y$ which respects the
structure maps, that is, $\theta$ intertwines the range and
source maps and satisfies
\[
(\id_\Lambda \times \theta) \circ \phi_X
= \phi_Y \circ (\theta \times \id_\Gamma).
\]
We call such a bijection $\theta$ an \emph{isomorphism} and
write $X \cong Y$; we denote the isomorphism class of a
$k$-morph $X$ by $[X]$.
\end{dfn}

We now introduce the notation associated with fibred products
of $k$-morphs, and then show in
Proposition~\ref{prp:composition} that the resulting object is
itself a $k$-morph.

\begin{ntn}\label{ntn:X*X}
Let $\Lambda_0$, $\Lambda_1$ and $\Lambda_2$ be $k$-graphs, and
let $X_i$ be a \morph{\Lambda_{i-1}}{\Lambda_i}\ with structure
maps $r_i, s_i$ and $\phi_i$ for $i = 1,2$. Let
\[
X_1 *_{\Lambda_1^0} X_2 = \{(x_1, x_2) \in X_1 \times X_2 : s(x_1)
= r(x_2)\}.
\]
Define $r : X_1 *_{\Lambda_1^0} X_2 \to \Lambda_0^0$ and $s :
X_1 *_{\Lambda_1^0} X_2 \to \Lambda_2^0$ by
\[
r(x_1, x_2) = r_1(x_1)\qquad\text{and}\qquad s(x_1, x_2) = s_2(x_2).
\]
To define $\phi : X_1 *_{\Lambda_1^0} X_2 *_{\Lambda_2^0}
\Lambda_2 \to \Lambda_0 *_{\Lambda_0^0} X_1 *_{\Lambda_1^0}
X_2$, fix $((x_1, x_2), \lambda_2) \in X_1 *_{\Lambda_1^0} X_2
*_{\Lambda_2^0} \Lambda_2$. Then $s_2(x_2) = r(\lambda_2)$, so
$\phi_2(x_2, \lambda_2) = (\lambda_1, x_2')$ for some
$\lambda_1 \in \Lambda_1$ and $x_2' \in X_2$. Moreover,
$r(\lambda_1) = r_2(x_2) = s_1(x_1)$, so $\phi_1(x_1,
\lambda_1) = (\lambda_0, x_1')$ for some $\lambda_0 \in
\Lambda_0$ and $x_1' \in X_1$ with $s_1(x_1') = r_2(x_2')$. We
define
\begin{equation}\label{eq:X*Yphi}
\phi((x_1, x_2), \lambda_2) = (\lambda_0, (x_1', x_2')).
\end{equation}
\end{ntn}

\begin{prop}\label{prp:composition}
With the notation above, $X_1 *_{\Lambda_1^0} X_2$ is a
\morph{\Lambda_0}{\Lambda_2}. Moreover, the isomorphism class
$[X_1 *_{\Lambda_1^0} X_2]$ depends only on the isomorphism
classes $[X_1]$ and $[X_2]$.
\end{prop}
\begin{proof}
Conditions \eqref{it:morph d}--\eqref{it:morph r} of
Definition~\ref{dfn:k-morphs} are easily checked
using~\eqref{eq:X*Yphi} and that $X_1$ and $X_2$ are
$k$-morphs, so we need only check~\eqref{it:phi assoc}.

Fix a composable pair $\mu_2, \nu_2 \in \Lambda_2$ and $(x_1,
x_2) \in X_1 *_{\Lambda_1^0} X_2$ such that $s(x_1, x_2) =
r(\mu_2)$. Let $x_i', x_i'' \in X_i$ for $i = 1,2$ and $\mu_i,
\nu_i \in \Lambda_i$ for $i = 0,1$ be the unique elements such
that
\begin{align}
\phi_i(x_i, \mu_i) &= (\mu_{i-1}, x_i') \label{eq:muthru} \\
\phi_i(x_i',\nu_i) &= (\nu_{i-1}, x_i''). \label{eq:nuthru}
\intertext{so that by~\eqref{eq:X*Yphi},}
\phi((x_1, x_2), \mu_2) &= (\mu_0, (x_1', x_2')) \label{eq:muresult}\\
\phi((x_1', x_2'), \nu_2) &= (\nu_0, (x_1'', x_2'')) \label{eq:nuresult}
\end{align}
By definition, $\phi((x_1, x_2),\mu_2\nu_2)$ is calculated as
follows: we write $\phi_2(x_2, \mu_2\nu_2) = (\lambda, y)$ and
then write $\phi_1(x_1, \lambda) = (\lambda', y')$; we then
have $\phi((x_1, x_2), \mu_2\nu_2) = (\lambda', (y', y))$. By
Definition~\ref{dfn:k-morphs}\eqref{it:phi assoc} for the
$k$-morph $X_2$, and equations
\eqref{eq:muthru}~and~\eqref{eq:nuthru}, we have $\lambda =
\mu_1\nu_1$ and $y = x_2''$. Now
Definition~\ref{dfn:k-morphs}\eqref{it:phi assoc} for the
$k$-morph $X_1$, and equations
\eqref{eq:muthru}~and~\eqref{eq:nuthru} force $\lambda' =
\mu_0\nu_0$ and $y' = x_1''$. That is, $\phi((x_1,
x_2),\mu_2\nu_2) = (\mu_0\nu_0, (x_1'', x_2''))$. Combining
this with equations \eqref{eq:muresult}~and~\eqref{eq:nuresult}
shows that $X_1 *_{\Lambda^0_1} X_2$ satisfies
Definition~\ref{dfn:k-morphs}\eqref{it:phi assoc}, and
therefore is a $k$-morph.

For the last statement, one checks that if $\theta_1 : X_1 \to
X_1'$ and $\theta_2 : X_2 \to X_2'$ are isomorphisms, then
$\theta_1 \times \theta_2$ is an isomorphism of $X_1
*_{\Lambda_1^0} X_2$ onto $X_1' *_{\Lambda_1^0} X_2'$.
\end{proof}

\begin{rmks}\label{rmks:composition}\ \par %\label{rmk:coveringcomposition}
\begin{enumerate}\renewcommand{\theenumi}{\roman{enumi}}
\item \label{it:coveringcomposition} Let $\Lambda_0$,
    $\Lambda_1$ and $\Lambda_2$ be $k$-graphs, and let $q :
    \Lambda_2 \to \Lambda_1$ and $p : \Lambda_1 \to
    \Lambda_0$ be coverings. Then $p \circ q$ is a covering
    of $\Lambda_0$ by $\Lambda_2$. Furthermore,
    $\sideset{_p}{}{\opX} *_{\Lambda_1^0}
    \sideset{_q}{}{\opX} = \{(q(v), v) : v \in
    \Lambda^0_2\}$. One can easily check that $\theta : v
    \mapsto (q(v),v)$ determines an isomorphism of
    $k$-morphs $ \sideset{_{p \circ q}}{}{\opX} \cong \sideset{_p}{}{\opX}
    *_{\Lambda_1^0} \sideset{_q}{}{\opX}$.
\item \label{it:pandq} Fix coverings $p : \Gamma \to
    \Lambda_1$ and $q : \Gamma \to \Lambda_2$. Let
    $\sideset{_p}{}{\opX}$, $X_q$ and
    $\sideset{_p}{_q}{\opX}$ be as in parts (ii),
    (iii)~and~(iv) respectively of Examples~\ref{ex:first}.
    Then we have $\sideset{_p}{_q}{\opX} \cong
    \sideset{_p}{}{\opX} *_{\Gamma^0}
    \sideset{}{_q}{\opX}$: the isomorphism $\theta$ is
    defined by $\theta(x) = (x, x)$.
\item \label{it:id} Let $\Lambda$ and $\Gamma$ be
    $k$-graphs, and let $X$ be a {\LGmph}. Then there are
    isomorphisms
    \[
    I_\Lambda *_{\Lambda^0} X \cong X \cong X *_{\Gamma^0} I_\Gamma.
    \]
    determined by $(r(x), x) \mapsto x$ and $(x, s(x))
    \mapsto x$.
\end{enumerate}
\end{rmks}

To state the next lemma, we describe the fibred product of $n$
$k$-morphs.

Let $\Lambda_0$, $\Lambda_1$, \dots, $\Lambda_n$ be $k$-graphs,
and let $X_i$ be a \morph{\Lambda_{i-1}}{\Lambda_i}\ for $i =
1, \dots, n$. Let
\[
X_1 *_{\Lambda_1^0} \cdots *_{\Lambda_{n-1}^0} X_n
 = \{(x_1, \cdots, x_n) \in X_1 \times \cdots \times X_n :
     s(x_{i-1}) = r(x_i)\text{ for } 1 \le i \le n\}.
\]
Define structure maps $r, s, \phi$ associated to $X = X_1
*_{\Lambda_1^0} \cdots *_{\Lambda_{n-1}^0} X_n$ as follows. We
set $r(x_1, \dots, x_n) = r(x_1)$ and $s(x_1, \dots, x_n) =
s(x_n)$. Given $((x_1, \dots, x_n), \lambda_n) \in X
*_{\Lambda_n^0} \Lambda_n$, let $\lambda_i \in \Lambda_i$ and
$x_i' \in X_i$ be the unique elements such that
\[
\phi_i(x_i, \lambda_i) = (\lambda_{i-1}, x_i')
    \text{ for $1 \le i  \le n$}.
\]
Since each $X_i$ is a $k$-morph, we have $s(x_i') =
s(\lambda_i) = r(x_{i+1}')$ for $1 \le i \le n-1$. So we define
$\phi$ by $\phi((x_1, \dots, x_n), \lambda_n) = (\lambda_0,
(x_1', \dots, x_n'))$.

\begin{lem}\label{lem:assoc}
With notation as above, $X = X_1 *_{\Lambda_1^0} \dots
*_{\Lambda_{n-1}^0} X_n$ is a \morph{\Lambda_0}{\Lambda_n}. For
$1 < m \le n$, there is an isomorphism
\[
(X_1 *_{\Lambda_1^0} \cdots *_{\Lambda_{m-2}^0} X_{m-1})
  *_{\Lambda_{m-1}^0}
 (X_m *_{\Lambda_m^0} \cdots *_{\Lambda_{n-1}^0} X_n)
\cong X
\]
implemented by $\theta_{m, n-m}((x_1, \dots, x_{m-1}),(x_m,
\dots, x_n)) = (x_1, \dots, x_n)$.
\end{lem}
\begin{proof}
%We first show
We first show that $X$ is a $k$-morph: properties
\eqref{it:morph d}--\eqref{it:morph r} of
Definition~\ref{dfn:k-morphs} are clear from the definition of
$\phi_{X}$, and Definition~\ref{dfn:k-morphs}\eqref{it:phi
assoc} is established by iterating an argument similar to that
of Proposition~\ref{prp:composition}.

It is easy to see that the bijection $\theta_{m, n-m}$
determines an isomorphism.
\end{proof}

\begin{ntn}\label{ntn:*-powers}
Let $X$ be a $\Lambda$ endomorph. For $n \ge 2$, we write
$X^{*n}$ for the $k$-morph
\[
X^{*n} := \overbrace{X *_{\Lambda^0} X *_{\Lambda^0}
\cdots *_{\Lambda^0} X}^{n\text{ terms}}.
\]
By $X^{*1}$, we mean $X$, and by $X^{*0}$, we mean $I_\Lambda$.
\end{ntn}

\begin{thm}\label{thm:category}
There is a category $\Mm_k$ such that: $\Obj(\Mm_k)$ is the
class of $k$-graphs; $\Hom_{\Mm_k}(\Gamma,\Lambda)$ is the set
of isomorphism classes of {\LGmph}s; the identity morphism
associated to $\Lambda \in \Obj(\Mm_k)$ is $[I_\Lambda]$; and
the composition map
\[
\Hom_{\Mm_k}(\Lambda_1, \Lambda_0) \times \Hom_{\Mm_k}(\Lambda_2, \Lambda_1) \to \Hom_{\Mm_k}(\Lambda_2,
\Lambda_0)
\]
is defined by $([X_1],[X_2]) \mapsto [X_1 *_{\Lambda_1^0}
X_2]$.
\end{thm}
\begin{proof}
Remark~\ref{rmks:composition}\eqref{it:id}, shows that the
$[I_\Lambda]$ act as identity morphisms.
Proposition~\ref{prp:composition} shows that the composition
map is well-defined. Since Lemma~\ref{lem:assoc} gives
\[
(X_1 *_{\Lambda_1^0} X_2) *_{\Lambda_2^0} X_3 \cong X \cong
X_1 *_{\Lambda_1^0} (X_2 *_{\Lambda_2^0} X_3)
\]
whenever the expressions make sense, the composition map is
also associative.
\end{proof}

\begin{rmk}\label{rmk:invertible}
In light of the preceding theorem, it is natural to ask
when the isomorphism class of a $k$-morph is an invertible
morphism of $\Mm_k$.

By an abuse of terminology, we will say that a \LGmph\ $X$ is
\emph{invertible} if $[X]$ is invertible in $\Mm_k$; that is,
if there is a \morph{\Gamma}{\Lambda}\ $Y$ such that $X
*_{\Gamma^0} Y \cong I_\Lambda$ and $Y *_{\Lambda^0} X \cong
I_\Gamma$. We claim that $X$ is invertible if and only if $X$
is isomorphic to $X(\alpha)$ for some isomorphism $\alpha :
\Gamma \to \Lambda$.

To see this, we first show that the range and source maps on
$X$ are bijections. Certainly the range and source maps on $X$
and $Y$ are surjective. In particular, $s : X \to \Gamma^0$ and
$r : Y \to \Gamma^0$ are surjective, and hence the projections
from $X *_{\Gamma^0} Y$ to $X$ and $Y$ are surjective. As $X
*_{\Gamma^0} Y$ is isomorphic to $I_\Lambda$, the range and
source maps on $X *_{\Gamma^0} Y$ are bijections. Since the
range map on $X *_{\Gamma^0} Y$ is defined by first projecting
onto the first coordinate in $X *_{\Gamma^0} Y$ and then
applying the range map on $X$, it follows that the range map on
$X$ is bijective. A similar argument shows that the source map
on $Y$ is bijective. Applying the same argument to $Y
*_{\Lambda^0} X \cong I_\Lambda$ shows that the range map on
$Y$ and the source map on $X$ are both bijective.

We may now define $\alpha : \Gamma \to \Lambda$ as follows.
Given $\gamma \in \Gamma$, there is a unique $x \in X$ with
$s(x) = r(\gamma)$. We then have $\phi(x,\gamma) = (\lambda,
x')$ for some $\lambda \in \Lambda$ and $x' \in X$, and we
define $\alpha(\gamma) = \lambda$. The properties of $\phi$ can
be used to show that $\alpha$ is an isomorphism of $k$-graphs,
and it is straightforward to check that $X \cong X(\alpha)$.
\end{rmk}

\section{Linking graphs}\label{sec:linking}

In this section, we define the notion of a linking graph
$\Sigma$ for a \LGmph\ $X$. Linking graphs generalise the
$(k+1)$-graphs $\Lambda \cs{p} \Gamma$ built from covering maps
$p$ in~\cite{KPS}.

We begin by showing how appropriate inclusions of $k$-graphs
$\Lambda$ and $\Gamma$ in a $(k+1)$-graph $\Sigma$ can be used
to manufacture a \LGmph\ $X$. This will provide a template for
linking graphs (see Definition~\ref{dfn:linking graph}).

\begin{ntn}\label{ntn:linking template}
Let $\Sigma$ be a $(k+1)$-graph, and let $\Lambda$ and $\Gamma$
be $k$-graphs. Let $\iota : \NN^k \to \NN^{k+1}$ be the
inclusion $\iota(n) = (n,0)$. Suppose that $i : \Lambda \sqcup
\Gamma \to \Sigma$ is an $\iota$-quasimorphism (where $\sqcup$
denotes a disjoint union) such that $i$ induces a $k$-graph
isomorphism $\Lambda \sqcup \Gamma \cong
\subgrph{\Sigma}{\iota}$. Suppose further that for all $\alpha
\in \Sigma^{e_{k+1}}$, $s(\alpha) \in i(\Gamma^0)$ and
$r(\alpha) \in i(\Lambda^0)$.

Let $X(\Lambda, \Gamma, \Sigma, i) = \Sigma^{e_{k+1}}$. Define
structure maps on $X = X(\Lambda, \Gamma, \Sigma, i)$ as
follows. For $\alpha \in X$, $r_X(\alpha)$ is the unique vertex
$v \in \Lambda^0$ such that $i(v) = r(\alpha)$, and similarly,
$s_X(\alpha)$ is the unique $w \in \Gamma^0$ satisfying $i(w) =
s(\alpha)$. For $(\alpha,\gamma) \in X *_{\Gamma^0} \Gamma$,
the factorisation property in $\Sigma$ ensures that $\alpha
i(\gamma) = i(\lambda) \alpha'$ for some unique $\lambda \in
\Lambda^{d(\gamma)}$ and $\alpha' \in X$, and we define
$\phi_X(\alpha, \gamma) = (\lambda, \alpha')$.
\end{ntn}

\begin{lem}\label{lem:linking template}
With the notation just established, $X(\Lambda, \Gamma, \Sigma,
i)$ is a \LGmph.
\end{lem}
\begin{proof}
Properties \eqref{it:morph d}--\eqref{it:morph r} of
Definition~\ref{dfn:k-morphs} are clear because $\Sigma$ is a
$(k+1)$-graph and $i$ is a quasimorphism.
Property~\eqref{it:phi assoc} follows from the associativity of
composition in $\Sigma$.
\end{proof}

\begin{dfn}\label{dfn:linking graph}
Let $\Lambda, \Gamma$ be $k$-graphs and let $X$ be a \LGmph.
Suppose that $\Sigma$, $\iota$ and $i$ are as in
Notation~\ref{ntn:linking template}. We say that the pair
$(\Sigma, i)$ is a \emph{linking graph} for $X$ if the
\LGmph\ $X(\Lambda, \Gamma, \Sigma, i)$ is isomorphic to $X$.
\end{dfn}

In practice, we will just say that $\Sigma$ is a linking graph
for $X$, leaving $i$ implicit.

\begin{example}
Let $p : \Gamma \to \Lambda$ be a covering of $k$-graphs, and
let $\sideset{_p}{}{\opX}$ be the $k$-morph described in
Examples~\ref{ex:first}\eqref{it:pX}. Then the $(k+1)$-graph
$\Lambda \cs{p} \Gamma$ of \cite[Proposition~2.6]{KPS} is a
linking graph for $\sideset{_p}{}{\opX}$.
\end{example}

Note that if $\Sigma$ is a linking graph for a $k$-morph, then
necessarily $d_\Sigma(\Sigma) \subset \NN^k \times \{0, 1\}$.

\begin{prop}\label{prp:linking graph exists}
Let $\Lambda, \Gamma$ be $k$-graphs and let $X$ be a \LGmph.
Then there exists a linking graph for $X$, and this linking
graph is unique up to isomorphism.
\end{prop}
\begin{proof}
As a set, we define $\Sigma = \Lambda \sqcup \Gamma \sqcup
(\Lambda *_{\Lambda^0} X)$. We endow $\Sigma$ with the
structure of a $(k+1)$-graph as follows.  First set $\Sigma^0 =
\Lambda^0 \sqcup \Gamma^0$. The restrictions of $r_\Sigma$ and
$s_\Sigma$ to $\Lambda \sqcup \Gamma$ are inherited from the
range and source maps on $\Lambda$ and $\Gamma$. For $\sigma
\in \Lambda \sqcup \Gamma$, set $d_\Sigma(\sigma) = (d(\sigma),
0)$. For $(\lambda, x) \in \Lambda *_{\Lambda^0} X$, let
$r_\Sigma(\lambda, x) = r(\lambda)$, $s_\Sigma(\lambda, x) =
s_X(x)$ and $d_\Sigma(\lambda, x) = (d(\lambda), 1)$.  Now fix
$\sigma_1, \sigma_2 \in \Sigma$ such that $s_\Sigma(\sigma_1) =
r_\Sigma(\sigma_2)$. We must define the composition $\sigma_1
\sigma_2$. There are three cases to consider.
\begin{enumerate}
\item If $\sigma_1, \sigma_2 \in \Gamma \sqcup \Lambda$,
    their composition as elements of $\Sigma$ is computed
    in $\Lambda \sqcup \Gamma$.
\item If $\sigma_1 = \mu \in \Lambda$, and $\sigma_2 =
    (\nu,x) \in \Lambda *_{\Lambda^0} X$, we define
    $\sigma_1\sigma_2 = (\mu\nu, x) \in \Lambda
    *_{\Lambda^0} X$.
\item If $\sigma_1 = (\mu, x) \in \Lambda *_{\Lambda^0} X$
    and $\sigma_2  = \nu \in \Gamma$, we write $\phi(x,
    \nu) = (\nu', x')$, and define $\sigma_1\sigma_2 =
    (\mu\nu', x') \in \Lambda *_{\Lambda^0} X$.
\end{enumerate}
Associativity follows from
Definition~\ref{dfn:k-morphs}\eqref{it:phi assoc}.

It is straightforward to check that $d_\Sigma : \Sigma \to
\NN^{k+1}$ is a functor. To show that $\Sigma$ has the
factorisation property, fix $\sigma \in \Sigma$, and suppose
$d_\Sigma(\sigma) = m + n$; we must show that there exist
unique paths $\tau \in \Sigma^m$ and $\rho \in \Sigma^n$ with
$\sigma = \tau\rho$. By definition of $d_\Sigma$, we have
$d_\Sigma(\sigma)_{k+1} \le 1$. If $m_{k+1} = n_{k+1} = 0$,
then $\sigma \in \Lambda \sqcup \Gamma$, and the factorisation
property in $\Lambda \sqcup \Gamma$ produces the desired paths
$\tau$ and $\rho$. If $m_{k+1} = 0$ and $n_{k+1} = 1$, then
$\sigma \in \Lambda *_{\Lambda^0} X$, say $\sigma = (\lambda,
x)$. By the factorisation property in $\Lambda$, there is a
unique factorisation $\lambda = \mu\nu$ where $d(\mu) = (m_1,
\dots, m_k)$, and then $\tau = \mu$ and $\rho = (\nu, x)$ are
the desired paths. Finally, suppose that $m_{k+1} = 1$ and
$n_{k+1} = 0$, and let $p \in \NN^k$ be the element such that
$m = (p, 1)$. Again, write $\sigma = (\lambda, x) \in \Lambda
*_{\Lambda^0} X$. Use the factorisation property in $\Lambda$
to write $\lambda = \mu\nu$ where $d(\mu) = p$. We have $(\nu,
x) \in \Lambda *_{\Lambda^0} X$, so $\phi(\nu,x) = (x', \nu')$
for some $x' \in X$ and $\nu' \in \Gamma$. One checks that
$\tau = (\mu,x')$ and $\rho = \nu'$ are the desired paths.

For uniqueness, let $\Sigma'$ be a linking graph for $X$. Then
there are a quasi-morphism $i' : \Lambda \sqcup \Gamma \to
\Sigma'$ satisfying the conditions set forth in
Notation~\ref{ntn:linking template} and an isomorphism $\psi :
X \to X(\Lambda, \Gamma, \Sigma', i')$. We define $\tilde\psi :
\Sigma \to \Sigma'$ as follows. For $\sigma \in \Lambda \sqcup
\Gamma \subset \Sigma$, we set $\tilde\psi(\sigma) =
i'(\sigma)$; and for $\sigma = (\lambda, x) \in \Lambda
*_{\Lambda^0} X$, we set $\tilde\psi(\lambda,x) =
i'(\lambda)\psi(x)$. One then checks that $\tilde\psi$ is an
isomorphism.
\end{proof}

\begin{rmk}\label{rmk:isomorphic linking graphs}
If $X$ and $X'$ are isomorphic $k$-morphs, then any linking
graph for $X$ is by definition also a linking graph for $X'$.
Hence Proposition~\ref{prp:linking graph exists} implies that,
up to isomorphism, there is a unique linking graph for each
isomorphism class of $k$-morphs.
\end{rmk}

\section{\texorpdfstring{$\Gamma$}{Gamma}-systems and \texorpdfstring{$\Gamma$}{Gamma}-bundles}\label{sec:systems}

In this section, we describe a generalisation, based on
$k$-morphs and linking graphs, of the $(k+1)$-graphs
$\tgrphlim(\Lambda_n; p_n)$ constructed in \cite{KPS} from a
sequence of coverings $p_n : \Lambda_{n+1} \to \Lambda_n$ of
$k$-graphs. The idea is that the sequence $(p_n)^\infty_{n=0}$
of coverings determines a consistent collection of $k$-morphs
\[
X_{(m,n)} = \sideset{_{p_{m}}}{}{\opX} *_{\Lambda^0_{m+1}} \dots *_{\Lambda_{n-1}^0} \sideset{_{p_{n-1}}}{}{\opX}
\]
indexed by pairs $(m,n) \in \NN \times \NN$ such that $m \le
n$. Recall from~\cite{KP} that such pairs are morphisms in the
$1$-graph $\Omega_1$ whose vertices are identified with $\NN$.
That is, a sequence of coverings as in \cite{KPS} gives rise to
a consistent collection of $k$-morphs indexed by a $1$-graph.

We generalise this situation by replacing $\Omega_1$ with an
$l$-graph $\Gamma$. Given a consistent collection (which we
call a $\Gamma$-system) of $k$-morphs indexed by paths in
$\Gamma$, we construct a $(k+l)$-graph $\Sigma$ which we call a
$\Gamma$-bundle. As we shall see, each $\Gamma$-system $X$
determines a functor $F_X$ from $\Gamma$ to $\Mm_k$. Na\"ively,
one might expect to be able to recover $\Sigma$ from $F_X$,
thus obviating the need to discuss $\Gamma$-systems at all. It
turns out, however, that this is not the case: not only do
there exist functors from which we may build two non-isomorphic
$(k+l)$-graphs, but there also exist functors from which no
$(k+l)$-graph may be built (see Examples~\ref{egs:not just
functors}).

\begin{dfn}\label{dfn:Gamma system}
Let $\Gamma$ be an $l$-graph, and let $k \ge 0$. Fix
\begin{itemize}
\item for each vertex $v \in \Gamma^0$ a $k$-graph
    $\Lambda_v$;
\item for each $\gamma \in \Gamma$ a
    \morph{\Lambda_{r(\gamma)}}{\Lambda_{s(\gamma)}}\
    $X_\gamma$; and
\item for each composable pair
    $\alpha,\beta$ in $\Gamma$ an isomorphism
    $\theta_{\alpha,\beta} : X_{\alpha}
    *_{\Lambda_{s(\alpha)}^0} X_{\beta} \to
    X_{\alpha\beta}$.
\end{itemize}
Suppose that $\Lambda$, $X$ and $\theta$ have the following
properties:
\begin{enumerate}
\item\label{it:Xid} for each $v \in \Gamma^0$, $X_v =
    I_{\Lambda_v}$,
\item\label{it:sys theta} for each $\gamma \in \Gamma$, the
    isomorphisms $\theta_{r(\gamma),\gamma}$ and
    $\theta_{\gamma, s(\gamma)}$ are those of
    Remark~\ref{rmks:composition}(\ref{it:id}), and
\item\label{it:sys assoc} for each composable triple
    $\alpha,\beta,\gamma \in \Gamma$, the following diagram
    commutes.
\[\begin{CD}
X_{\alpha} *_{\Lambda_{s(\alpha)}^0} X_{\beta} *_{\Lambda_{s(\beta)}^0} X_{\gamma}
 @>\theta_{\alpha,\beta} \times \id_{X_{\gamma}}>>
X_{\alpha\beta} *_{\Lambda_{s(\beta)}^0} X_{\gamma} \\
 @V\id_{X_{\alpha}} \times \theta_{\beta,\gamma}VV   @V\theta_{\alpha\beta,\gamma}VV \\
X_{\alpha} *_{\Lambda_{s(\alpha)}^0} X_{\beta\gamma}
 @>\hspace{1em}\theta_{\alpha, \beta\gamma}\hspace{1em}>>
X_{\alpha\beta\gamma}
\end{CD}\]
\end{enumerate}
Then we say that \emph{$X$ is a $\Gamma$-system of $k$-morphs
with data $\Lambda, \theta$} or just that $X$ is a
$\Gamma$-system, in which case $\Lambda$ and $\theta$ are
implicit.
\end{dfn}

\begin{rmk}\label{rmk:systems gives functor}
Given a $\Gamma$-system $X$ of $k$-morphs, there is a functor
$F_X$ from $\Gamma$ to $\Mm_k$ determined by $F_X(\gamma) =
[X_\gamma]$; in particular, the object map satisfies $F^0_X(v)
= \Lambda_v$. However, the $\Gamma$-system $X$ contains more
information than the functor $F_X$: the $\Gamma$-system picks
out a concrete representative $X_\gamma$ of each isomorphism
class $F_X(\gamma)$ and a compatible system of concrete
isomorphisms $X_{\alpha} *_{\Lambda^0_{s(\alpha)}} X_\beta
\cong X_{\alpha\beta}$ implementing the compositions
$F_X(\alpha)F_X(\beta) = F_X(\alpha\beta)$. We show that this
distinction is important in Examples~\ref{egs:not just
functors}.
\end{rmk}

\begin{examples} \label{ex:systems}\ \hfill \\[-2.5ex]
\begin{enumerate}\renewcommand{\theenumi}{\roman{enumi}}
\item\label{it:1-edge sys} Fix $k$-graphs $\Lambda$ and
    $\Gamma$ and a \LGmph\ $X$. Let $E$ be the unique
    $1$-graph with a single edge $e$ and two vertices
    $r(e)$ and $s(e)$. Setting $\Lambda_{r(e)} = \Lambda$,
    $\Lambda_{s(e)} = \Gamma$, $X_{r(e)} = I_\Lambda$,
    $X_{s(e)} = I_\Gamma$, and $X_e = X$, we obtain an
    $E$-system of $k$-morphs.
    \setcounter{mycount}{\value{enumi}}
\item \label{it:endo system} Fix a $k$-graph $\Lambda$ and
    a $\Lambda$ endomorph $X$. Recall from
    \cite[Examples~1.7(iii)]{KP} that $T_1$ denotes the
    unique $1$-graph with a single edge $e$ and a single
    vertex $v = r(e) = s(e)$. The degree functor is an
    isomorphism of $T_1$ onto $\NN$, and we use it to
    identify the two. Let $\Lambda_{0} = \Lambda$, and for
    $n \ge 0$ let $X_{n} = X^{*n}$ as in
    Notation~\ref{ntn:*-powers}. Then the isomorphisms
    $\theta_{m,n} : X_m *_{\Lambda^0} X_n \to X_{m+n}$ of
    Lemma~\ref{lem:assoc} give this collection of
    $k$-morphs the structure of a $T_1$-system.
\item\label{it:skew prod sys} Fix an $l$-graph $\Gamma$, a
    countable discrete group $G$, and a functor $c : \Gamma
    \to G$. For each $v \in \Gamma^0$, let $\Lambda_v$ be
    the $0$-graph such that $\Lambda^0_v = G$. For $\gamma
    \in \Gamma$, there is an automorphism
    $\alpha_{c(\gamma)} : \Lambda_{s(\gamma)} \to
    \Lambda_{r(\gamma)}$ determined by
    $\alpha_{c(\gamma)}(g) = c(\gamma)g$ for $g \in
    \Lambda^0_{s(\gamma)}$. For each $\gamma \in \Gamma$,
    let $X_\gamma$ be the $0$-morph $X(\alpha_{c(\gamma)})$
    arising from $\alpha_{c(\gamma)}$ as in
    Example~\ref{ex:first}(\ref{it:X(alpha)}). The
    isomorphisms $X(\alpha_{c(\gamma\gamma')}) \cong
    X(\alpha_{c(\gamma)}) *_{G} X(\alpha_{c(\gamma')})$
    described in
    Remark~\ref{rmks:composition}(\ref{it:coveringcomposition}),
    give this collection of $0$-morphs the structure of a
    $\Gamma$-system, denoted $X(c)$, of $0$-morphs.
\item Let $\Cov_k$ denote the category whose objects are
    $k$-graphs and whose morphisms are $k$-graph coverings,
    and let $\Gamma$ be an $l$-graph. Let $F$ be a functor
    from $\Gamma$ to $\Cov_k$. For each $v \in \Gamma^0$,
    let $\Lambda_v$ be the $k$-graph $F^0(v)$ and for each
    $\gamma \in \Gamma$, let $X_\gamma =
    \sideset{_{F(\gamma)}}{}{\opX}$. Then the isomorphisms
    $\theta_{\gamma,\gamma'} : X_\gamma
    *_{\Lambda^0_{s(\gamma)}} X_{\gamma'} \to
    X_{\gamma\gamma'}$ of
    Remark~\ref{rmks:composition}(\ref{it:coveringcomposition})
    give this collection of $k$-morphs the structure of a
    $\Gamma$-system. For instance, the covering systems of
    \cite{KPS} give rise to $\Omega_1$-systems.
\item\label{it:PRRS} Let $\Sigma$ be a $2$-graph satisfying
    the hypotheses of \cite[Theorem~3.1]{PRRS}. That is,
    each vertex lies on a unique simple cycle in the graph
    whose edges are $\Sigma^{e_2}$, and the graph with
    edge-set $\Sigma^{e_1}$ contains no cycles. The simple
    cycles in edges in $\Sigma^{e_2}$ determine an
    equivalence relation on vertices in $\Sigma$ by $v \sim
    w$ if and only if $v\Sigma^{ne_2} w \not= \emptyset$
    for some $n \in \NN$. We write $[v]$ for the
    equivalence class of $v \in \Sigma^0$ under this
    relation. There is a $1$-graph $\Gamma$ with
    \[
    \Gamma^0 = \{[v] : v \in \Sigma^0\}\quad\text{ and }\quad
    \Gamma^1 = \{([v],[w]) : [v],[w] \in \Sigma^0/\sim,
                             [v]\Sigma^{e_1}[w] \not=\emptyset\}
    \]
    where $r([v],[w]) = [v]$ and $s([v],[w]) = [w]$. For
    $[v] \in \Gamma^0$, we denote by $\Lambda_{[v]}$ the
    sub-1-graph of $\Sigma$ such that $\Lambda_{[v]}^1 =
    [v]\Sigma^{e_2}[v]$; so each $\Lambda_{[v]}$ is the
    path-category of a simple cycle. For each path $\gamma
    \in \Gamma$, we define $X_\gamma = r(\gamma)
    \Sigma^{d(\gamma)e_1} s(\gamma)$, and we endow it with
    the range and source maps $r : X_\gamma \to
    \Lambda^0_{r(\gamma)}$ and $s : X_\gamma \to
    \Lambda^0_{s(\gamma)}$ inherited from $\Sigma$, and
    with the map $\phi : X_\gamma *_{s(\gamma)}
    \Lambda_{s(\gamma)} \to \Lambda_{r(\gamma)}
    *_{r(\gamma)} X_\gamma$ determined by the factorisation
    property in $\Sigma$. Then each $X_\gamma$ is a
    \morph{\Lambda_{r(\gamma)}}{\Lambda_{s(\gamma)}},
    composition in $\Sigma$ defines isomorphisms
    $\theta_{\alpha,\beta} : X_{\alpha} *_{s(\alpha)}
    X_{\beta} \to X_{\alpha\beta}$ and from this structure
    we obtain a $\Gamma$-system of $1$-morphs which we
    denote by $X(\Sigma)$. If $\Sigma$ is a rank-2 Bratteli
    diagram as in \cite[Section~4]{PRRS}, then $\Gamma$ is
    the path-category of a Bratteli diagram.
\end{enumerate}
\end{examples}

The next step is to show how to assemble a $(k+l)$-graph from
the data contained in a $\Gamma$-system of $k$-graphs. This
construction should simultaneously generalise the linking
graphs of the previous section, and the construction of the
skew-product of a $k$-graph by a group.

The model for this construction is the following prototypical
$\Gamma$-system which generalises the construction of
Example~\ref{ex:systems}(\ref{it:PRRS}).

\begin{ntn}\label{ntn:prototype Gamma system}
Let $\Sigma$ be a $(k+l)$-graph, and let $\Gamma$ be an
$l$-graph. Let $\pi : \NN^{k+l} \to \NN^l$ denote the
projection onto the last $l$ coordinates; that is $\pi(m,n) =
n$. Suppose $f : \Sigma \to \Gamma$ is a $\pi$-quasimorphism
which restricts to a surjection of $\Sigma^0$ onto $\Gamma^0$.
Let $\imath : \NN^k \to \NN^{(k+l)}$ be
the inclusion $\imath(m) = (m,0)$ and let $\jmath : \NN^l \to
\NN^{k+l}$ be the inclusion $\jmath(n) = (0,n)$. For each $v
\in \Gamma^0$, we define $\Lambda_v = f^{-1}(v)$ which is a
sub-$k$-graph of $\Sigma^{\imath}$. For $\gamma \in \Gamma$, we
define $X(f)_\gamma = f^{-1}(\gamma) \cap \Sigma^\jmath$ (note
that for each $x \in X(f)_\gamma$, $d_\Sigma(x) = (0,
d(\gamma))$). Each $X(f)_\gamma$ becomes a
\morph{\Lambda_{r(\gamma)}}{\Lambda_{s(\gamma)}}\ under the
range, source and factorisation maps inherited from $\Sigma$.
Composition in $\Sigma$ defines maps
\[
\theta_{\gamma_1, \gamma_2} : X(f)_{\gamma_1}
*_{\Lambda^0_{s(\gamma_1)}} X(f)_{\gamma_2} \to
X(f)_{\gamma_1\gamma_2}
\]
for each composable pair $(\gamma_1, \gamma_2)$ in $\Gamma$.
Moreover, under these structure maps, $X(f)$ becomes a
$\Gamma$-system of $k$-graphs: conditions
(\ref{it:Xid})~and~(\ref{it:sys theta}) of
Definition~\ref{dfn:Gamma system} are satisfied by definition,
and condition~(\ref{it:sys assoc}) follows from associativity
of composition in $\Sigma$.
\end{ntn}

We will show that every $\Gamma$-system is isomorphic to one of
the form $X(f)$ for some $\pi$-quasimorphism $f : \Sigma \to
\Gamma$. We must first make clear what we mean by an
isomorphism of $\Gamma$-systems.

\begin{dfn}\label{dfn:Gammasysiso}
Let $\Gamma$ be a $l$-graph. Suppose that $(\Lambda, X,
\theta)$ and $(\Xi, Y, \psi)$ are $\Gamma$-systems of
$k$-graphs. An isomorphism from $X$ to $Y$ consists of
$k$-graph isomorphisms $h^0_v : \Lambda_v \to \Xi_v$ and
bijections $h_\gamma : X_\gamma \to Y_\gamma$ which intertwine
all the structure maps. That is:
\begin{enumerate}
\item for each $\gamma \in \Gamma$, $s \circ h_\gamma =
    h^0_{s(\gamma)} \circ s$, $r \circ h_\gamma =
    h^0_{r(\gamma)} \circ r$ and
    \[
        (h^0_{r(\gamma)} \times h_\gamma) \circ \phi_{X_\gamma}
         = \phi_{Y_\gamma} \circ (h_\gamma \times h^0_{s(\gamma)});\text{ and}
    \]
\item for every composable pair $\alpha,\beta \in \Gamma$,
    \[
        h_{\alpha\beta} \circ \theta_{\alpha,\beta}
            = \psi_{\alpha,\beta} \circ (h_\alpha \times h_\beta).
    \]
\end{enumerate}
\end{dfn}

We can now say what we mean by a $\Gamma$-bundle for a
$\Gamma$-system.

\begin{dfn}
Let $\Gamma$ be an $l$-graph, and let $X$ be a $\Gamma$-system
of $k$-graphs. Let $\pi : \NN^{k+l} \to \NN^l$ be the
projection onto the last $l$ coordinates. A
\emph{$\Gamma$-bundle} for $X$ is a $(k+l)$-graph $\Sigma$
endowed with a $\pi$-quasimorphism $f : \Sigma \to \Gamma$
which restricts to a surjection of $\Sigma^0$ onto $\Gamma^0$
such that the $\Gamma$-system $X(f)$ of
Notation~\ref{ntn:prototype Gamma system} is isomorphic to $X$.
We call the $\pi$-quasimorphism $f$ the \emph{bundle map} for
the $\Gamma$-system.
\end{dfn}

\begin{rmk}\label{rmk:the magic of the definition}
Our formulation of $\Gamma$-bundles emphasises that the
construction of Notation~\ref{ntn:prototype Gamma system} is
prototypical: given a $\pi$-quasimorphism $f : \Sigma \to
\Gamma$, the $(k+l)$-graph $\Sigma$ is automatically a
$\Gamma$-bundle (with bundle map $f$) for the resulting
$\Gamma$-system $X(f)$.
\end{rmk}

\begin{rmk}\label{rmk:skewgraphstructure}
Let $X$ be a $\Gamma$-system of $k$-graphs, and suppose that
$\Sigma$ is a $\Gamma$-bundle for $X$. The maps $h^0_v$ and
$h_\gamma$ of Definition~\ref{dfn:Gammasysiso} determine
injective quasimorphisms $h^0_v : \Lambda_v \to \Sigma$ for
each $v \in \Gamma^0$, and inclusions $h_\gamma : X_\gamma \to
\Sigma^{\jmath(d(\gamma))}$. Moreover, the factorisation
property implies that every element of $\sigma \in \Sigma$ can
be expressed as $\sigma = h_{f(\sigma)}(x)
h^0_{s(f(\sigma))}(\lambda)$ for unique elements, $x \in
X_{f(\sigma)}$ and $\lambda \in \Lambda_{s(f(\sigma))}$.
\end{rmk}

We now show that every $\Gamma$-system admits a
$\Gamma$-bundle.

\begin{thm}\label{thm:skewgraphexists}
Let $\Gamma$ be an $l$-graph and let $X$ be a $\Gamma$-system
of $k$-morphs. Then there exists a $\Gamma$-bundle $\Sigma$ for
$X$. Moreover, $\Sigma$ is unique up to isomorphism: that is,
if $\Psi$ is another $\Gamma$-bundle for $X$, then there is an
isomorphism $\Sigma \cong \Psi$ of $(k+l)$-graphs which
intertwines the bundle maps on $\Sigma$ and $\Psi$.
\end{thm}
\begin{proof}
Throughout this proof, for $\gamma \in \Gamma$, we will write
$\phi_\gamma$ for the isomorphism $\phi_{X_\gamma} : X_{\gamma}
*_{\Lambda_{s(\gamma)}^0} \Lambda_{s(\gamma)} \to \Lambda_{r(\gamma)}
*_{\Lambda_{r(\gamma)}^0} X_{\gamma}$ associated with the $k$-morph
$X_\gamma$.

We must first construct a $\Gamma$-bundle for $X$, and then
show that it is unique. We begin by constructing the
$(k+l)$-graph $\Sigma$. Define a set $\Sigma$ by
\[
\Sigma = \{(\lambda,\gamma,x) : \gamma \in \Gamma,
\lambda \in \Lambda_{r(\gamma)}, x \in X_\gamma\}.
\]
Define $d : \Sigma \to \NN^{k+l}$ by $d(\lambda,\gamma,x) =
(d(\lambda), d(\gamma))$. We write $\Sigma^p$ for the set
$d^{-1}(p) \subset \Sigma$ for each $p \in \NN^{k+l}$. By
condition~(\ref{it:Xid}) of Definition~\ref{dfn:Gamma system},
for $v \in \Gamma^0$, $X_v = I(\Lambda_v)$ is equal as a set to
$\Lambda_v^0$. Hence $\Sigma^0 = \{(u,v,u) : v \in \Gamma^0, u
\in \Lambda_{v}^0\}$. To simplify notation and to help
distinguish vertices from arbitrary paths, we will discard the
redundant $u$, and write $(u,v)$ for the element $(u,v,u)$ of
$\Sigma^0$.

Define $r,s : \Sigma \to \Sigma^0$ by
\[
r(\lambda,\gamma,x) = (r(\lambda), r(\gamma))\quad\text{and}\quad
s(\lambda,\gamma,x) = (s(x), s(\gamma)).
\]
Suppose that $s(\lambda_0, \gamma_0, x_0) = r(\lambda_1, \gamma_1,
x_1)$, and define $\lambda_1'$ and $x'_0$ by  $\phi_{\gamma_0}(x_0, \lambda_1) =
(\lambda'_1, x'_0)$. We define
\[
(\lambda_0, \gamma_0, x_0)(\lambda_1, \gamma_1, x_1)
 = (\lambda_0\lambda'_1, \gamma_0\gamma_1, \theta_{\gamma_0,\gamma_1}(x'_0,x_1)).
\]
It is easy to check that the triple on the right lies in
$\Sigma$ and that
\[
d(\lambda_0\lambda'_1, \gamma_0\gamma_1, \theta_{\gamma_0,\gamma_1}(x'_0,x_1)) =
d(\lambda_0, \gamma_0, x_0) + d(\lambda_1, \gamma_1, x_1).
\]

We aim to show that $\Sigma$ is a $(k+l)$-graph when endowed
with these structure maps. We must check that the composition
we have defined is associative, and that under this
composition, $(\Sigma,d)$ satisfies the factorisation property.
To see that the composition is associative, we fix a composable
triple $(\lambda_0, \gamma_0, x_0), (\lambda_1, \gamma_1, x_1),
(\lambda_2, \gamma_2, x_2)$ of elements of $\Sigma$. Let $v_i =
r(\gamma_i)$ and $v_{i+1} = s(\gamma_i)$ for $i = 0,1,2$.
Define $x_0', x_1', x_1'', \lambda_1', \lambda_2'$ and
$\lambda_2''$ by
\begin{align*}
\phi_{\gamma_0}(x_0, \lambda_1) &= (\lambda'_1, x'_0),\\
\phi_{\gamma_1}(x_1, \lambda_2) &= (\lambda'_2, x'_1),\quad\text{and}\\
\phi_{\gamma_0}(x_0', \lambda_2') &= (\lambda_2'', x_0'').
\end{align*}
We may visualise the situation as follows:
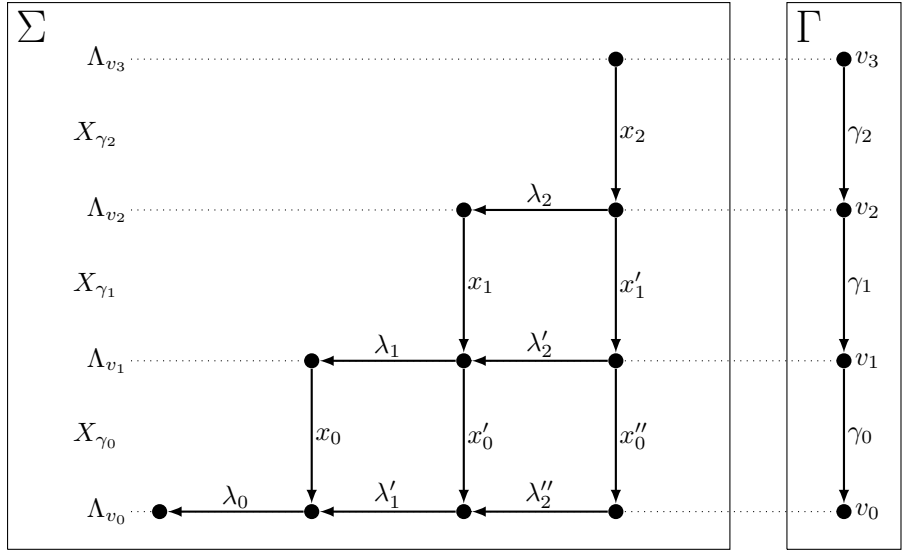
\[
\begin{tikzpicture}
    \draw[style=thin] (8.25,-0.5)--(9.75,-0.5)--(9.75,6.75)--(8.25,6.75)--(8.25,-0.5);
    \node[inner sep=3pt, anchor=north west] at (8.25,6.75) {\Large$\Gamma$};
    \draw[style=thin] (-2,-0.5)--(7.5,-0.5)--(7.5,6.75)--(-2,6.75)--(-2,-0.5);
    \node[inner sep=3pt, anchor=north west] at (-2,6.75) {\Large$\Sigma$};
    \draw[style=thin, style=dotted] (9,0)--(-.4,0);
    \draw[style=thin, style=dotted] (9,2)--(-.4,2);
    \draw[style=thin, style=dotted] (9,4)--(-.4,4);
    \draw[style=thin, style=dotted] (9,6)--(-.4,6);
    \node[inner sep=1pt, anchor=east] at (-0.5,1) {\footnotesize$X_{\gamma_0}$};
    \node[inner sep=1pt, anchor=east] at (-0.5,3) {\footnotesize$X_{\gamma_1}$};
    \node[inner sep=1pt, anchor=east] at (-0.5,5) {\footnotesize$X_{\gamma_2}$};
    \node[inner sep=1pt, anchor=west] at (-1,0) {\footnotesize$\Lambda_{v_0}$};
    \node[inner sep=1pt, anchor=west] at (-1,2) {\footnotesize$\Lambda_{v_1}$};
    \node[inner sep=1pt, anchor=west] at (-1,4) {\footnotesize$\Lambda_{v_2}$};
    \node[inner sep=1pt, anchor=west] at (-1,6) {\footnotesize$\Lambda_{v_3}$};
    \node[inner sep=2pt, circle, fill=black] (00) at (0,0) {};
    \node[inner sep=2pt, circle, fill=black] (10) at (2,0) {};
    \node[inner sep=2pt, circle, fill=black] (20) at (4,0) {};
    \node[inner sep=2pt, circle, fill=black] (30) at (6,0) {};
    \node[inner sep=2pt, circle, fill=black] (11) at (2,2) {};
    \node[inner sep=2pt, circle, fill=black] (21) at (4,2) {};
    \node[inner sep=2pt, circle, fill=black] (31) at (6,2) {};
    \node[inner sep=2pt, circle, fill=black] (22) at (4,4) {};
    \node[inner sep=2pt, circle, fill=black] (32) at (6,4) {};
    \node[inner sep=2pt, circle, fill=black] (33) at (6,6) {};
    \node[inner sep=2pt, circle, fill=black] (G0) at (9,0) {};
    \node[inner sep=2pt, circle, fill=black] (G1) at (9,2) {};
    \node[inner sep=2pt, circle, fill=black] (G2) at (9,4) {};
    \node[inner sep=2pt, circle, fill=black] (G3) at (9,6) {};
    \draw[style=thick, -latex] (G3.south)--(G2.north) node[pos=0.5, anchor=west,inner sep=1pt] {\footnotesize$\gamma_2$};
    \draw[style=thick, -latex] (G2.south)--(G1.north) node[pos=0.5, anchor=west,inner sep=1pt] {\footnotesize$\gamma_1$};
    \draw[style=thick, -latex] (G1.south)--(G0.north) node[pos=0.5, anchor=west,inner sep=1pt] {\footnotesize$\gamma_0$};
    \node[inner sep=1pt, anchor=west] at (G3.east) {\footnotesize$v_3$};
    \node[inner sep=1pt, anchor=west] at (G2.east) {\footnotesize$v_2$};
    \node[inner sep=1pt, anchor=west] at (G1.east) {\footnotesize$v_1$};
    \node[inner sep=1pt, anchor=west] at (G0.east) {\footnotesize$v_0$};
    \draw[style=thick, -latex] (10.west)--(00.east) node[pos=0.5, anchor=south,inner sep=1pt] {\footnotesize$\lambda_0$};
    \draw[style=thick, -latex] (20.west)--(10.east) node[pos=0.5, anchor=south,inner sep=1pt] {\footnotesize$\lambda_1'$};
    \draw[style=thick, -latex] (30.west)--(20.east) node[pos=0.5, anchor=south,inner sep=1pt] {\footnotesize$\lambda_2''$};
    \draw[style=thick, -latex] (21.west)--(11.east) node[pos=0.5, anchor=south,inner sep=1pt] {\footnotesize$\lambda_1$};
    \draw[style=thick, -latex] (31.west)--(21.east) node[pos=0.5, anchor=south,inner sep=1pt] {\footnotesize$\lambda_2'$};
    \draw[style=thick, -latex] (32.west)--(22.east) node[pos=0.5, anchor=south,inner sep=1pt] {\footnotesize$\lambda_2$};
    \draw[style=thick, -latex] (33.south)--(32.north) node[pos=0.5, anchor=west,inner sep=1pt] {\footnotesize$x_2$};
    \draw[style=thick, -latex] (32.south)--(31.north) node[pos=0.5, anchor=west,inner sep=1pt] {\footnotesize$x_1'$};
    \draw[style=thick, -latex] (31.south)--(30.north) node[pos=0.5, anchor=west,inner sep=1pt] {\footnotesize$x_0''$};
    \draw[style=thick, -latex] (22.south)--(21.north) node[pos=0.5, anchor=west,inner sep=1pt] {\footnotesize$x_1$};
    \draw[style=thick, -latex] (21.south)--(20.north) node[pos=0.5, anchor=west,inner sep=1pt] {\footnotesize$x_0'$};
    \draw[style=thick, -latex] (11.south)--(10.north) node[pos=0.5, anchor=west,inner sep=1pt] {\footnotesize$x_0$};
\end{tikzpicture}
\]
To prove associativity, we must show that the products
\[
\big((\lambda_0, \gamma_0, x_0)
(\lambda_1, \gamma_1, x_1)\big) (\lambda_2, \gamma_2, x_2)
\quad\text{ and }\quad
(\lambda_0, \gamma_0, x_0) \big((\lambda_1, \gamma_1, x_1) (\lambda_2,
\gamma_2, x_2)\big)
\]
coincide. We begin by calculating the first of these. First,
notice that the pair $((x'_0, x_1), \lambda_2)$ belongs to
$(X_{\gamma_0} *_{\Lambda^0_{v_1}} X_{\gamma_1})
*_{\Lambda^0_{v_2}} \Lambda_{v_2}$. We have
\begin{align}
\phi_{X_{\gamma_0} *_{\Lambda^0_{v_1}} X_{\gamma_1}}((x'_0, x_1), \lambda_2)
 &= (\lambda_2'', (x_0'',x_1')) \label{eq:identity1}
\intertext{by definition of $\phi_{X_{\gamma_0} *_{\Lambda^0_{v_1}} X_{\gamma_1}}$
and of the elements $\lambda_2''$, $x_0''$ and $x_1'$ given above. Since
$\theta_{\gamma_0, \gamma_1} : X_{\gamma_0} *_{\Lambda^0_{v_1}} X_{\gamma_1}
\to X_{\gamma_0\gamma_1}$ is an isomorphism of $k$-morphs,}
\phi_{\gamma_0\gamma_1}(\theta_{\gamma_0,\gamma_1}(x_0',x_1),\lambda_2)
 &= (\id_{\Lambda_{v_0}} \times \theta_{\gamma_0,\gamma_1})
    (\phi_{X_{\gamma_0} *_{\Lambda^0_{v_1}} X_{\gamma_1}}((x'_0, x_1), \lambda_2)) \nonumber\\
 &= (\lambda_2'', \theta_{\gamma_0,\gamma_1}(x_0'',x_1')) \label{eq:identity2}
\intertext{by~\eqref{eq:identity1}. Therefore,}
\big((\lambda_0, \gamma_0, x_0) (\lambda_1, \gamma_1, x_1)\big) (\lambda_2, \gamma_2, x_2)
 &= (\lambda_0\lambda'_1, \gamma_0\gamma_0, \theta_{\gamma_0,\gamma_1}(x_0',x_1))(\lambda_2, \gamma_2, x_2)\nonumber\\
 &= ((\lambda_0\lambda'_1)\lambda''_2,\,
     (\gamma_0\gamma_1)\gamma_2,\,
     \theta_{\gamma_0\gamma_1, \gamma_2}(\theta_{\gamma_0,\gamma_1}(x_0'',x_1'),x_2),\nonumber
\intertext{where the second step uses~\eqref{eq:identity2}. Similar calculations show that}
(\lambda_0, \gamma_0, x_0) \big((\lambda_1, \gamma_1, x_1) (\lambda_2, \gamma_2, x_2)\big)
 &= (\lambda_0(\lambda'_1\lambda''_2),
     \gamma_0(\gamma_1\gamma_2),
     \theta_{\gamma_0, \gamma_1\gamma_2}(x_0'', \theta_{\gamma_1,\gamma_2}(x_1', x_2)).\nonumber
\end{align}
Associativity in $\Sigma$ now follows from associativity in
$\Lambda_{v_0}$ and $\Gamma$, and property~(\ref{it:sys assoc})
of Definition~\ref{dfn:Gamma system}.

To establish the factorisation property in $\Sigma$, fix $m,p \in
\NN^k$ and $n,q \in \NN^l$ and an element $(\lambda,\gamma,x) \in
\Sigma^{(m+p, n+q)}$. By the factorisation properties in
$\Lambda_{r(\gamma)}$ and in $\Gamma$, there are unique
factorisations $\lambda = \lambda_0\lambda_1$ and $\gamma =
\gamma_0\gamma_1$ where $d(\lambda_0) = m$, $d(\lambda_1) = p$,
$d(\gamma_0) = n$ and $d(\gamma_1) = q$. Since $\theta_{\gamma_0,
\gamma_1}$ is an isomorphism, there are unique elements $x_0 \in
X_{\gamma_0}$ and $x_1 \in X_{\gamma_1}$ such that $x =
\theta_{\gamma_0,\gamma_1}(x_0, x_1)$. As $\phi_{\gamma_0}$ is also a
bijection, there are unique elements $\lambda_1' \in
\Lambda_{s(\gamma_0)}$ and $x_0' \in X_{\gamma_0}$ such that
$\phi_{\gamma_0}(x_0', \lambda_1') = (\lambda_1,x_0)$. We then have
$d(\lambda_0, \gamma_0, x_0') = (m,p)$ and $d(\lambda_1', \gamma_1,
x_1) = (n,q)$, and
\[
(\lambda,\gamma,x) = (\lambda_0, \gamma_0, x_0')(\lambda_1', \gamma_1, x_1)
\]
by definition. Uniqueness is clear. We have now established
that $\Sigma$ is a $(k+l)$-graph.

The formula $f(\lambda,\gamma,x) = \gamma$ defines a
$\pi$-quasimorphism from $\Sigma$ onto $\Gamma$. Let $\imath :
\NN^k \to \NN^{k+l}$ and $\jmath : \NN^l \to \NN^{k+l}$ be as
in Notation~\ref{ntn:prototype Gamma system}. For each $v \in
\Gamma^0$, $f^{-1}(v) = \{(\lambda, v, s(\lambda)) : \lambda
\in \Lambda_v\}$, and $h^0_v : f^{-1}(v) \to \Lambda_v$ defined
by $h_v^0(\lambda, v, s(\lambda)) = \lambda$ is an isomorphism
of $k$-graphs. For each $\gamma \in \Gamma$, $f^{-1}(\gamma)
\cap \Sigma^\jmath = \{(r(x), \gamma, x) : x \in X_\gamma\}$,
and $h_\gamma : f^{-1}(\gamma) \cap \Sigma^\jmath \to X_\gamma$
defined by $(r(x), \gamma, x) = x$ is a bijection. Routine
calculations show that these maps satisfy conditions
(1)~and~(2) of Definition~\ref{dfn:Gammasysiso}. Hence $\Sigma$
is a $\Gamma$-bundle for $X$ as claimed when equipped with the
bundle map $f : \Sigma \to \Gamma$.

It remains to establish the uniqueness of the $\Gamma$-bundle.
Suppose that $\Psi$ is another $\Gamma$-bundle for $X$ with
bundle map $g$. So we have isomorphisms $h^0_v : \Lambda_v \to
g^{-1}(v)$ for each $v \in \Gamma^0$ and bijections $h_\gamma :
X_\gamma \to g^{-1}(\gamma) \cap \Psi^\jmath$ determining an
isomorphism of $\Gamma$-systems. Define $H : \Sigma \to \Psi$
by
\[
H(\lambda,\gamma,x) =
h^0_{r(\gamma)}(\lambda)h_\gamma(x)\quad\text{ for all $(\lambda,x,\gamma)
\in \Sigma$.}
\]
The factorisation property in $\Psi$ ensures that each $\psi
\in \Psi$ with $d(\psi) = (m,n)$ can be written uniquely as
$\psi = \psi_m\psi_n$ where $d(\psi_m) = (m,0)$ and $d(\psi_n)
= (0,n)$. We then have $\psi_m \in g^{-1}(g(r(\psi)))$, and
$\psi_n \in g^{-1}(g(\psi)) \cap \Psi^\jmath$. Bijectivity of
$H$ follows from this. It is clear that $H$ respects the degree
map, and intertwines the range and source maps. A
straightforward calculation shows that it also respects
composition, and hence is an isomorphism of $(k+l)$-graphs.
\end{proof}

\begin{rmk}
Theorem~\ref{thm:skewgraphexists} implies that it makes sense
to talk about \emph{the} $\Gamma$-bundle for a $\Gamma$-system
$X$, and we shall frequently do so. Unless specified otherwise,
the bundle is denoted $\Sigma$ and the bundle map is denoted
$f$.
\end{rmk}

\begin{examples} \label{ex:skew-products}\ \hfill \\[-2.5ex]
\begin{enumerate}\renewcommand{\theenumi}{\roman{enumi}}
\item Let $\Lambda$ and $\Gamma$ be $k$-graphs, and let $X$
    be a \LGmph. As noted in
    Example~\ref{ex:systems}(\ref{it:1-edge sys}), this
    corresponds to an $E$-system of $k$-graphs where $E$ is
    the $1$-graph with a single edge $e$ and two vertices
    $r(e)$ and $s(e)$. An $E$-bundle for this $E$-system
    amounts to a linking graph for $X$.
\item\label{it:endoskewprod} Let $\Lambda$ be a $k$-graph,
    and $X$ a $\Lambda$ endomorph. As noted in
    Example~\ref{ex:systems}(\ref{it:endo system}), this
    corresponds to a $T_1$-system of $k$-graphs. We shall
    denote the $T_1$-bundle for this system by $\Lambda
    \times_X \NN$. We call $\Lambda \times_X \NN$ the
    \emph{endomorph skew-graph for $X$}. Every
    $(k+1)$-graph arises this way: given a $(k+1)$-graph
    $\Sigma$, with $\Lambda = \Sigma^\iota$ and $X =
    \Sigma^{e_{k+1}}$ as in Example~\ref{ex:first}(vi),
    $\Sigma$ is isomorphic to $\Lambda \times_X \NN$.
\item\label{it:alphacrossprod} Let $\Lambda$ be a
    $k$-graph, and let $\alpha$ be an automorphism of
    $\Lambda$. Let $X = X(\alpha)$ be the associated
    $\Lambda$ endomorph. In this case, the endomorph
    skew-graph $\Lambda \times_X \NN$ discussed in the
    preceding example is the same as the crossed-product
    $(k+1)$-graph $\Lambda \times_\alpha \ZZ$ of
    \cite{FPS}. More generally, let $T_l$ denote the
    $l$-graph isomorphic to $\NN^l$, and suppose that
    $\alpha$ is an action of $\ZZ^l$ by automorphisms of
    $\Lambda$. Let $X_n = X(\alpha_n)$ for each $n \in
    T_l$, and let $\theta_{m,n} : X_m *_{\Lambda^0} X_n \to
    X_{m+n}$ be the isomorphism of
    Remark~\ref{rmks:composition}(\ref{it:coveringcomposition}).
    Then the $X_n$ form a $T_l$-system $X(\alpha)$, and the
    crossed product $(k+l)$-graph $\Lambda \times_\alpha
    \ZZ^l$ described in \cite{FPS} is a $T_l$-bundle for
    $X(\alpha)$.
\item Fix a functor $c$ from an $l$-graph $\Gamma$ to a
    group $G$, and construct from this a $\Gamma$-system
    $X(c)$ of $0$-morphs as in
    Example~\ref{ex:systems}(\ref{it:skew prod sys}). Then
    the skew-product $k$-graph $G \times_c \Gamma$
    of~\cite{KP} is a $\Gamma$-bundle for the
    $\Gamma$-system $X(c)$; the bundle map is the functor
    $f(g, \gamma) := \gamma$.
\item Let $\Sigma$ be a $2$-graph satisfying the hypotheses
    of \cite[Theorem~3.1]{PRRS}, and let $\Gamma$ and $X$
    be as in Example~\ref{ex:systems}(\ref{it:PRRS}).
%%     If $F : \NN^2 \to \NN^2$ denotes the homomorphism
%%     $F(m,n) = (n,m)$, then a $\Gamma$-bundle for $X$ is
%%     $F$-quasimorphic to $\Lambda$.
    Then $\Sigma$ is a $\Gamma$-bundle for $X$ and the
    bundle map is the natural quotient map $\Lambda \to
    \Sigma$. In particular, we may regard a rank-2 Bratteli
    diagram $\Sigma$ as a bundle of cycle-graphs over the
    path-category of a conventional Bratteli diagram.
\end{enumerate}
\end{examples}

Our $\Gamma$-bundle construction is quite general: the next
proposition shows that every $k$-graph $\Lambda$ is a
$T_k$-bundle for some $T_k$-system. It is a strong point of our
formulation of $\Gamma$-bundles that the proof of this result
is almost trivial (see Remark~\ref{rmk:the magic of the
definition}).

\begin{prop}
Let $\Lambda$ be a $k$-graph.
Let $T_k$ denote the $k$-graph isomorphic to $\NN^k$.
Then the degree map on $\Lambda$ determines an
$\id_k$-quasimorphism (also denoted $d$) from $\Lambda$ onto
$T_k$. In particular $\Lambda$ is a $T_k$-bundle (with bundle
map $d$) for the $T_k$-system $X(d)$ of $0$-morphs described in
Notation~\ref{ntn:prototype Gamma system}.
\end{prop}

\begin{rmk}
If $p : \Gamma \to \Lambda$ is a covering of $k$-graphs then
$p$ is an $\id_k$-quasimorphism. Hence $p$ induces a
$\Lambda$-system of $0$-morphs as in
Notation~\ref{ntn:prototype Gamma system}. Moreover, these
$0$-morphs are all invertible (see
Remark~\ref{rmk:invertible}).

Conversely, given a $k$-graph $\Lambda$, and a $\Lambda$-system of
invertible $0$-morphs, the bundle map associated to a
$\Lambda$-bundle for the system is a covering map.
\end{rmk}

We conclude the section by investigating the relationship
between $\Gamma$-systems and functors from $\Gamma$ into
$\Mm_k$. If $\Gamma$ is a $1$-graph, the two are essentially
the same thing.

\begin{prop}\label{prp:1-graph systems}
Let $\Gamma$ be a $1$-graph, and let $F : \Gamma \to \Mm_k$ be
a functor. Then there is, up to isomorphism, exactly one
$\Gamma$-system $X$ of $k$-graphs such that $F(\gamma) =
[X_\gamma]$ for all $\gamma \in \Gamma$.
\end{prop}
\begin{proof}
For each $v \in \Gamma^0$, set $\Lambda_{v} = F^0(v)$, and
$X_{v} = I(\Lambda_v)$. For each edge $e \in \Gamma^1$, fix a
$k$-morph $X_e$ such that $[X_e] = F(e)$. For $n \ge 2$ and a
path $\alpha = e_1 \cdots e_{n} \in \Gamma^n$, let $X_\alpha =
X_{e_1} *_{\Lambda^0_{s(e_1)}} X_{e_2} *_{\Lambda^0_{s(e_2)}}
\cdots *_{\Lambda^0_{s(e_{n-1})}} X_{e_n}$ as in
Lemma~\ref{lem:assoc}, and for composable $\alpha,\beta$, let
$\theta_{\alpha,\beta}$ be the isomorphism described in the
same lemma. It is easy to verify that this data determines a
$\Gamma$-system of $k$-morphs which induces the functor $F$.

Now suppose that $Y$ is another $\Gamma$ system of $k$-morphs
(with data $\Lambda$, $\psi$) which induces $F$. Let
$\psi_{\alpha,\beta} : Y_\alpha *_{\Lambda_{s(\alpha)}^0}
Y_\beta \to Y_{\alpha\beta}$ denote the isomorphisms in the
$\Gamma$-system $Y$. In particular, each $Y_\alpha$ is a
\morph{\Lambda_{r(\alpha)}}{\Lambda_{s(\alpha)}} which is
isomorphic to $X_\alpha$. For $v \in \Gamma^0$, let $h^0_v$
denote the identity map on $\Lambda_v$. For each $e \in
\Gamma^1$, we may fix an isomorphism $h_e : X_e \to Y_e$. By
induction on $n$, for $\alpha \in \Gamma^n$ and $f \in
\Gamma^1$ with $s(\alpha) = r(f)$, we may define an isomorphism
$h_{\alpha f}$ from $X_{\alpha f}$ to $Y_{\alpha f}$ by
\[
X_{\alpha f} \cong X_\alpha *_{\Lambda_{s(\alpha)}^0} X_f
\xrightarrow{h_\alpha \times h_f}
Y_\alpha *_{\Lambda_{s(\alpha)}^0} Y_f
\xrightarrow{\psi_{\alpha, h_f}}
Y_{\alpha f}.
\]
As the isomorphisms $h_\alpha$ are defined using the structure
maps in $X$ and $Y$, it is easy to check that they determine an
isomorphism of $\Gamma$-systems.
\end{proof}

\begin{examples}\label{egs:not just functors}
We cannot expect an analogue of Proposition~\ref{prp:1-graph systems}
to hold if $\Gamma$ is an $l$-graph with $l > 1$, as the following
two examples show.

\begin{enumerate} \renewcommand{\theenumi}{\roman{enumi}}
\item \label{it:non-unique system} Let $\Gamma = T_2$ (the
    $2$-graph isomorphic to $\NN^2$), and let $\Lambda$ be
    the $0$-graph consisting of a single vertex $v$. Each
    finite set is a $\Lambda$ endomorph when endowed with
    the only possible structure maps. In particular, each
    multiplicative map $x : \NN^2 \to \NN$ determines a
    functor from $T_2$ to $\Mm_0$: the image of $n \in
    \NN^2$ is (the isomorphism class of) the $\Lambda$
    endomorph $X_n$ with $x(n)$ elements.

    Example~6.1 of \cite{KP} describes two non-isomorphic
    $2$-graphs $\Lambda$ and $\Lambda'$ each with a single
    vertex, two edges of degree $e_1$ and two edges of
    degree $e_2$. As in Notation~\ref{ntn:prototype Gamma
    system}, $\Lambda$ and $\Lambda'$ determine
    $T_2$-systems $X$ and $X'$ such that $\Lambda$ is a
    $T_2$-bundle for $X$ and $\Lambda'$ is a $T_2$-bundle
    for $X'$. By Theorem~\ref{thm:skewgraphexists}, $X$ and
    $X'$ are non-isomorphic. However, $X$ and $X'$
    determine the same functor $F : T_2 \to \Mm_0$ with
    $F^0(0) = \Lambda$, namely the one corresponding to the
    multiplicative map $x : \NN^2 \to \NN$ given by $x(n) =
    2^{n_1 + n_2}$.

\item \label{it:Jacks} The following example is due to Jack
    Spielberg \cite{JS example}. We thank Jack for allowing
    us to reproduce it here. The following diagram
    represents a 3-coloured graph where the edges have
    colours $c_1$, $c_2$ and $c_3$; we draw $c_1$-coloured
    edges as solid lines, $c_2$-coloured edges as dashed
    lines, and $c_3$-coloured edges as dotted lines.
    
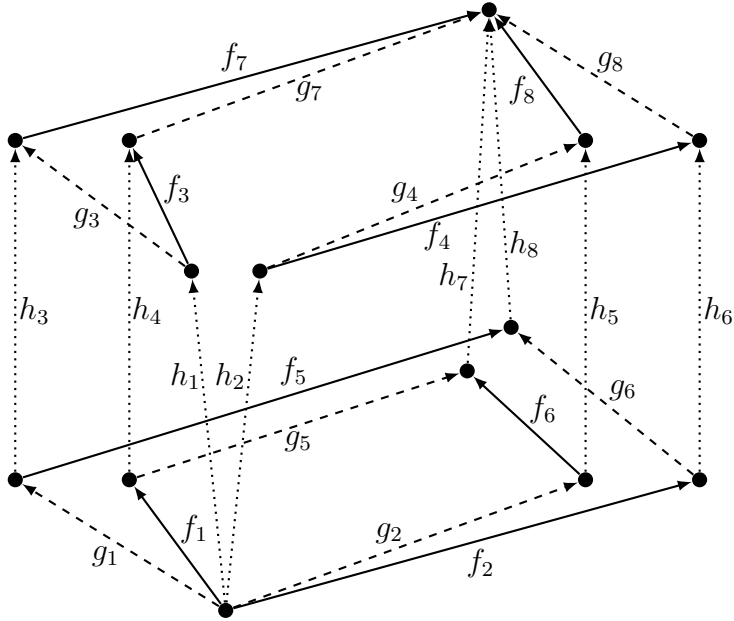
\[\begin{tikzpicture}[scale=1.5]
        \node[circle, fill=black, inner sep=2pt] (v1) at (0,0,0) {};%
        \node[circle, fill=black, inner sep=2pt] (v2) at (-3,0,-3) {};%
        \node[circle, fill=black, inner sep=2pt] (v3) at (-2,0,-3) {};%
        \node[circle, fill=black, inner sep=2pt] (v4) at (2,0,-3) {};%
        \node[circle, fill=black, inner sep=2pt] (v5) at (3,0,-3) {};%
        \node[circle, fill=black, inner sep=2pt] (v6a) at (0,0,-6.5) {};%
        \node[circle, fill=black, inner sep=2pt] (v6b) at (0,0,-5.5) {};%
        \node[circle, fill=black, inner sep=2pt] (w1a) at (-0.3,3,0) {};%
        \node[circle, fill=black, inner sep=2pt] (w1b) at (0.3,3,0) {};%
        \node[circle, fill=black, inner sep=2pt] (w2) at (-3,3,-3) {};%
        \node[circle, fill=black, inner sep=2pt] (w3) at (-2,3,-3) {};%
        \node[circle, fill=black, inner sep=2pt] (w4) at (2,3,-3) {};%
        \node[circle, fill=black, inner sep=2pt] (w5) at (3,3,-3) {};%
        \node[circle, fill=black, inner sep=2pt] (w6) at (0,3,-6) {};%
        \draw[style=thick, -latex] (v1)--(v3) node[pos=0.5, inner sep=0.5pt, anchor=south west] {$f_1$};%
        \draw[style=thick, -latex] (w1a)--(w3) node[pos=0.5, inner sep=0.5pt, anchor=south west] {$f_3$};%
        \draw[style=thick, -latex] (v2)--(v6a) node[pos=0.6, inner sep=1pt, anchor=south east] {$f_5$};%
        \draw[style=thick, -latex] (w2)--(w6) node[pos=0.5, inner sep=1pt, anchor=south east] {$f_7$};%
        \draw[style=thick, -latex] (v1)--(v5) node[pos=0.5, inner sep=1pt, anchor=north west] {$f_2$};%
        \draw[style=thick, -latex] (w1b)--(w5) node[pos=0.4, inner sep=1pt, anchor=north] {$f_4$};%
        \draw[style=thick, -latex] (v4)--(v6b) node[pos=0.5, inner sep=0.5pt, anchor=south west] {$f_6$};%
        \draw[style=thick, -latex] (w4)--(w6) node[pos=0.5, inner sep=0pt, anchor=north east] {$f_8$};%
        \draw[style=dashed, style=thick, -latex] (v1)--(v2) node[pos=0.5, inner sep=0pt, anchor=north east] {$g_1$};%
        \draw[style=dashed, style=thick, -latex] (w1a)--(w2) node[pos=0.5, inner sep=0.5pt, anchor=north east] {$g_3$};%
        \draw[style=dashed, style=thick, -latex] (v3)--(v6b) node[pos=0.5, inner sep=2pt, anchor=north] {$g_5$};%
        \draw[style=dashed, style=thick, -latex] (w3)--(w6) node[pos=0.5, inner sep=2pt, anchor=north] {$g_7$};%
        \draw[style=dashed, style=thick, -latex] (v1)--(v4) node[pos=0.5, inner sep=0.5pt, anchor=south east] {$g_2$};%
        \draw[style=dashed, style=thick, -latex] (w1b)--(w4) node[pos=0.5, inner sep=1pt, anchor=south east] {$g_4$};%
        \draw[style=dashed, style=thick, -latex] (v5)--(v6a) node[pos=0.5, inner sep=1pt, anchor=south west] {$g_6$};%
        \draw[style=dashed, style=thick, -latex] (w5)--(w6) node[pos=0.5, inner sep=0.5pt, anchor=south west] {$g_8$};%
        \draw[style=dotted, style=thick, -latex] (v1)--(w1a) node[pos=0.7, inner sep=0pt, anchor=east] {$h_1$};%
        \draw[style=dotted, style=thick, -latex] (v1)--(w1b) node[pos=0.7, inner sep=1pt, anchor=east] {$h_2$};%
        \draw[style=dotted, style=thick, -latex] (v2)--(w2) node[pos=0.5, inner sep=1pt, anchor=west] {$h_3$};%
        \draw[style=dotted, style=thick, -latex] (v3)--(w3) node[pos=0.5, inner sep=0.5pt, anchor=west] {$h_4$};%
        \draw[style=dotted, style=thick, -latex] (v4)--(w4) node[pos=0.5, inner sep=1pt, anchor=west] {$h_5$};%
        \draw[style=dotted, style=thick, -latex] (v5)--(w5) node[pos=0.5, inner sep=1pt, anchor=west] {$h_6$};%
        \draw[style=dotted, style=thick, -latex] (v6b)--(w6) node[pos=0.25, inner sep=1pt, anchor=east] {$h_7$};%
        \draw[style=dotted, style=thick, -latex] (v6a)--(w6) node[pos=0.25, inner sep=1pt, anchor=west] {$h_8$};%
    \end{tikzpicture}\]
    For distinct $1 \le i,j \le 3$, there is a unique
    range- and source-preserving bijection $\theta_{i,j}$
    between $c_ic_j$-coloured paths and $c_jc_i$-coloured
    paths. For example, the only $c_2c_1$-coloured path
    with the same range and source as the $c_1c_2$-coloured
    path $f_5g_1$ is $g_6f_2$, so $\theta_{1,2}(f_5g_1) =
    g_6f_2$. Thus, the factorisation rules in any $3$-graph
    with the skeleton pictured above must be implemented by
    the $\theta_{i,j}$. To see that no such 3-graph exists,
    we consider the two possible ways of reversing the
    colouring of the path $h_8g_6f_2$ using the
    $\theta_{i,j}$:
    \begin{align*}
    h_8g_6f_2 &\to h_8f_5g_1 \to f_7h_3g_1 \to f_7g_3h_1 \quad\text{ and}\\
    h_8g_6f_2 &\to g_8h_6f_2 \to g_8f_4h_2 \to f_8g_4h_2.
    \end{align*}
    Since $f_8 \not= f_7$, $g_3 \not= g_4$ and $h_1 \not=
    h_2$, the $\theta_{i,j}$ do not specify a valid
    collection of factorisation rules (see
    \cite[Section~2]{RSY1}).

    Let $\Gamma = T_3$ (the $3$-graph isomorphic to
    $\NN^3$), and let $\Lambda$ be the $0$-graph whose
    vertices are those in the diagram above. The sets $X_1
    := \{f_1, \dots f_8\}$, $X_2 := \{g_1, \dots, g_8\}$
    and $X_3 = \{h_1, \dots, h_8\}$ are
    $\Lambda$-endomorphs when endowed with the obvious
    structure maps. For $i \not= j$, $\theta_{i,j}$
    determines an isomorphism $X_i *_{\Lambda^0} X_j \cong
    X_j *_{\Lambda^0} X_i$, so there is a unique functor $F
    : \Gamma \cong \NN^3 \to \Mm_0$ such that $F^0(v) =
    \Lambda$ and $F(e_i) = [X_i]$ for $i = 1,2,3$. However,
    this functor is not determined by any $T_3$-system of
    $0$-morphs: the $T_3$-bundle for such a system would be
    a $3$-graph whose skeleton was the 3-coloured graph we
    started with.
    \end{enumerate}
\end{examples}

\section{\texorpdfstring{$C^*$}{C*}-correspondences and functoriality}\label{sec:functoriality}

\setcounter{equation}{-1}

In this section we consider how the constructions of the
preceding sections behave with respect to higher-rank graph
$C^*$-algebras. To keep the length of the paper down, we
restrict attention to \morph{\Lambda}{\Gamma}s such that
\begin{equation}\label{eq:regmorph}
\parbox{0.9\textwidth}{$\Lambda$ and $\Gamma$ are row-finite
$k$-graphs with no sources, $s : X \to \Gamma^0$ is surjective,
and $r : X \to \Lambda^0$ is surjective and finite-to-one.}
\tag{{\footnotesize\regsymbol}}
\end{equation}\stepcounter{equation}
This simplifying assumption ensures that the $\Gamma$-bundles
we construct are covered by the results of \cite{RSY1}.

To each $k$-morph $X$ satisfying~\malteseref, we associate a
$C^*$-correspondence $\Hh(X)$. The germ of this construction,
at least for $k$-morphs arising from covering maps, is present
in the proof of \cite[Proposition~3.2]{KPS}. However, here we
make it explicit and extend it to arbitrary $k$-morphs.

In Theorem~\ref{thm:functoriality}, we show that the assignment
$[X] \mapsto [\Hh(X)]$ determines a contravariant functor to
the category $\Cc$ whose objects are $C^*$-algebras and whose
morphisms are isomorphism classes of $C^*$-correspondences.
Theorem~\ref{thm:CP algebra} shows that when $X$ is a $\Lambda$
endomorph satisfying~\malteseref{}, the Cuntz-Pimsner algebra
of $\Hh(X)$ is isomorphic to the $(k+1)$-graph $C^*$-algebra
$C^*(\Lambda \times_X \NN)$.

\begin{prop}
There is a subcategory $\rfMm{k}$ of $\Mm_k$ whose objects are
row-finite $k$-graphs with no sources, and whose morphisms are
isomorphism classes of $k$-morphs $X$ satisfying~\malteseref{}.
\end{prop}
\begin{proof}
This follows from the observation that if $k$-morphs $X_1$ and
$X_2$ satisfy~\malteseref{}, then so does their product.
\end{proof}

\begin{lem}\label{lem:inclusions}
Let $\Lambda$ and $\Gamma$ be $k$-graphs and let $X$ be a \LGmph\
satisfying~\malteseref{}. Let $(\Sigma,i)$ be a linking graph
for $X$. Then $\Sigma$ is row-finite and locally convex, and
there are injective homomorphisms
\begin{align*}
i^*_\Lambda : C^*(\Lambda) \to C^*(\Sigma) &\text{ such that } i^*_\Lambda(s_\lambda) = s_{i(\lambda)},\text{ and} \\
i^*_\Gamma : C^*(\Gamma) \to C^*(\Sigma) &\text{ such that } i^*_\Gamma(s_\gamma) = s_{i(\gamma)}.
\end{align*}
The series $\sum_{v \in \Lambda^0} s_{i(v)}$ and $\sum_{w \in
\Gamma^0} s_{i(w)}$ converge strictly to complementary full
projections $P_\Lambda$ and $P_\Gamma$ in $MC^*(\Sigma)$. The
homomorphism $i^*_\Gamma$ induces an isomorphism $C^*(\Gamma)
\cong P_\Gamma C^*(\Sigma) P_\Gamma$. The homomorphism
$i^*_\Lambda$ induces an embedding $C^*(\Lambda)
\hookrightarrow P_\Lambda C^*(\Sigma) P_\Lambda$ which takes an
approximate identity for $C^*(\Lambda)$ to an approximate
identity for $P_\Lambda C^*(\Sigma) P_\Lambda$.
\end{lem}
\begin{proof}
The linking graph $\Sigma$ is row-finite because~\malteseref{}
ensures that $\Lambda$ and $\Gamma$ are both row-finite and the
range map on $X$ is finite-to-one. To see that $\Sigma$ is
locally convex, suppose that $e, f \in \Sigma$ satisfy $r(e) =
r(f)$, $d(e) = e_i$ and $d(f) = e_j$ where $1 \le i < j \le
k+1$. If $j \le k$, then $s(e)\Sigma^{e_j}$ and
$s(f)\Sigma^{e_i}$ are nonempty because $\Lambda$ and $\Gamma$
have no sources. If $d(f) = e_{k+1}$, then
$s(e)\Sigma^{e_{k+1}}$ is nonempty because $r : X \to
\Lambda^0$ is surjective by~\malteseref{}, and
$s(f)\Sigma^{e_i}$ is nonempty because $\Gamma$ has no sources.

The existence of homomorphisms $i^*_\Lambda$ and $i^*_\Gamma$
satisfying the required formulae follows from the universal
properties of $C^*(\Lambda)$ and $C^*(\Gamma)$, and their
injectivity follows from the gauge-invariant uniqueness theorem
\cite[Theorem~3.4]{KP}.

A standard argument (see for example
\cite[Proposition~3.2]{KPS}) shows that $P_\Lambda$ and
$P_\Gamma$ make sense and are complementary projections. To see
that $P_\Lambda$ is full, we fix a generator $s_\sigma$ of
$C^*(\Sigma)$ and show that $s_\sigma \in C^*(\Sigma) P_\Lambda
C^*(\Sigma)$. If $r(\sigma) \in i_\Lambda(\Lambda^0)$, then
$s_\sigma = P_\Lambda s_\sigma$; and if $r(\sigma) \in
i_\Gamma(\Gamma^0)$, then since the source map on $X$ is
surjective by~\malteseref{}, we have $s_\sigma = s^*_\alpha
P_\Lambda s_\alpha s_\sigma$ for some $\alpha \in
\Sigma^{e_{k+1}}$. To see that $P_\Gamma$ is full, we fix a
generator $s_\sigma$ of $C^*(\Sigma)$ and show that $s_\sigma
\in C^*(\Sigma) P_\Gamma C^*(\Sigma)$. If $s(\sigma) \in
i_\Gamma(\Gamma^0)$, then $s_\sigma = s_\sigma P_\Gamma$; and
if $s(\sigma) \in i_\Lambda(\Lambda^0)$, then since the range
map on $X$ is surjective and finite-to-one by~\malteseref{}, we
have
\[
s_\sigma = \sum_{\alpha \in s(\sigma)\Sigma^{e_{k+1}}} s_\sigma s_\alpha P_\Gamma s^*_\alpha.
\]

We have $P_\Gamma C^*(\Sigma) P_\Gamma = \clsp\{s_\alpha
s^*_\beta : \alpha,\beta \in \Sigma, r(\alpha),r(\beta) \in
i(\Gamma^0)\} = \clsp\{i^*_\Gamma(s_{\gamma_1}s_{\gamma_2}^*) :
\gamma_i \in \Gamma, s(\gamma_1) = s(\gamma_2)\}$, and it
follows from an application of the gauge-invariant uniqueness
theorem \cite[Theorem~4.1]{RSY1} that $i^*_\Gamma$ implements
the desired isomorphism $C^*(\Gamma) \cong P_\Gamma C^*(\Sigma)
P_\Gamma$. Since $\Sigma$ is row-finite and locally convex, the
Cuntz-Krieger relations in $C^*(\Sigma)$ ensure that the
partial isometries $\{s_{i_\Lambda(\lambda)} : \lambda \in
\Lambda\}$ form a Cuntz-Krieger $\Lambda$-family. Another
application of \cite[Theorem~4.1]{RSY1} then shows that
$i^*_\Lambda : C^*(\Lambda) \to P_\Lambda C^*(\Sigma)
P_\Lambda$ is injective. For each finite subset $F \subset
\Lambda^0$, let $p_F$ denote the projection $\sum_{v \in F}
s_v$. Then the net $(p_F)_{F \subset \Lambda^0\text{ finite}}$
is an approximate identity for $C^*(\Lambda)$, and
$(i^*_\Lambda(p_F))_{F \subset \Lambda^0\text{ finite}}$
converges strictly to $P_\Lambda$ by definition of $P_\Lambda$.
\end{proof}

\begin{dfn}\label{dfn:H(X)}
Resume the hypotheses of Lemma~\ref{lem:inclusions}. Let
$\Hh(X)$ denote the vector space $P_\Lambda C^*(\Sigma)
P_\Gamma$. Define a left action of $C^*(\Lambda)$ and a right
action of $C^*(\Gamma)$ on $\Hh(X)$ by
\[
a \cdot \xi \cdot b = i^*_\Lambda(a)\,\xi\,i^*_\Gamma(b)
\quad\text{ for $a \in C^*(\Lambda)$, $b \in C^*(\Gamma)$ and $\xi \in \Hh(X)$}
\]
where the product is taken in $C^*(\Sigma)$. Define $\langle
\cdot, \cdot \rangle_{C^*(\Gamma)} : \Hh(X) \times \Hh(X) \to
C^*(\Gamma)$ as follows: $\langle \xi,
\eta\rangle_{C^*(\Gamma)}$ is the unique element of
$C^*(\Gamma)$ such that
\[
\xi^*\eta = i^*_\Gamma(\langle \xi, \eta\rangle_{C^*(\Gamma)})
\]
where $\xi^*\eta$ is calculated in $C^*(\Sigma)$.
\end{dfn}

\begin{prop}\label{prp:H(X)}
Let $\Lambda$ and $\Gamma$ be $k$-graphs and let $X$ be a \LGmph\
satisfying~\malteseref{}. The space $\Hh(X)$ defined above
satisfies
\begin{equation}\label{eq:spanningHX}
\Hh(X) = \clsp\{s_{x i(\alpha)}s^*_{i(\beta)} : x \in X, \alpha,\beta \in \Gamma, s(x) = r(\alpha), s(\alpha) = s(\beta)\}.
\end{equation}
Under the operations defined above, $\Hh(X)$ is a full
nondegenerate $C^*(\Lambda)$--$C^*(\Gamma)$
$C^*$-correspondence, and the left-action is implemented by an
injective homomorphism of $C^*(\Lambda)$ into $\Kk(\Hh(X))$.
Moreover, the isomorphism class of $\Hh(X)$ depends only on the
isomorphism class of $X$.
\end{prop}

\begin{rmk}
If we do not insist on~\malteseref{}, but assume only that
$\Sigma$ is finitely aligned (so that $C^*(\Sigma)$ makes
sense), then Definition~\ref{dfn:H(X)} still specifies a
%% $(P_\Lambda C^*(\Sigma) P_\Lambda)$--$(P_\Gamma C^*(\Sigma)
%% P_\Gamma)$
$C^*(\Lambda)$--$C^*(\Gamma)$ correspondence $\Hh(X)$
satisfying~\eqref{eq:spanningHX}. However our proofs of the
remaining assertions of Proposition~\ref{prp:H(X)} and of
Theorems \ref{thm:functoriality}~and~\ref{thm:CP algebra} all
rely on~\malteseref{} via their dependence on
Lemma~\ref{lem:inclusions}.
\end{rmk}

\begin{proof}[Proof of Proposition~\ref{prp:H(X)}]
Fix a nonzero spanning element $s_\mu s^*_\nu$ of
$C^*(\Sigma)$. Then $s_\mu s^*_\nu \in P_\Lambda C^*(\Sigma)
P_\Gamma$ only if $r(\mu) \in i(\Lambda^0)$, $r(\nu) \in
i(\Gamma^0)$, and $s(\mu) = s(\nu)$. Since $r(\nu) \in
i(\Gamma^0)$ implies $s(\nu) \in \Gamma^0$, we have $\nu \in
i(\Gamma)$. Since $s(\mu) = s(\nu)$, we also have $\mu \in
i(\Lambda^0) \Sigma i(\Gamma^0)$, and the factorisation
property in $\Sigma$ forces $\mu = x i(\alpha)$ for some $x \in
X$ and $\alpha \in \Gamma$ with $s(x) = r(\alpha)$. This
establishes~\eqref{eq:spanningHX}.

Since $P_\Lambda$ and $P_\Gamma$ are complementary full
projections in $MC^*(\Sigma)$, Theorem~3.19 of~\cite{TFB}
implies that $\Hh(X)$ is an imprimitivity bimodule for the two
corners $P_\Lambda C^*(\Sigma) P_\Lambda \cong \Kk(\Hh(X))$ and
$P_\Gamma C^*(\Sigma) P_\Gamma \cong C^*(\Gamma)$. That
$\Hh(X)$ is full follows from the definition of an
imprimitivity bimodule (see \cite[Definition~3.1]{TFB}). The
injective homomorphism $C^*(\Lambda) \to\Kk(\Hh(X))$ comes from
the embedding $C^*(\Lambda) \hookrightarrow P_\Lambda
C^*(\Sigma) P_\Lambda$ induced by $i^*_\Lambda$ and the
identification $P_\Lambda C^*(\Sigma) P_\Lambda \cong
\Kk(\Hh(X))$.  Since $i^*_\Lambda$ maps an approximate identity
for $C^*(\Lambda)$ to an approximate identity for $P_\Lambda
C^*(\Sigma) P_\Lambda$, $\Hh(X)$ is nondegenerate.

The final statement follows from Remark~\ref{rmk:isomorphic
linking graphs}.
\end{proof}

For the following, recall from
Section~\ref{sec:correspondences} that $\Cc$ denotes the
category whose objects are $C^*$-algebras and whose morphisms
are isomorphism classes of $C^*$-correspondences.

\begin{thm}\label{thm:functoriality}
For each $k \ge 0$, the assignments $\Lambda \mapsto
C^*(\Lambda)$ and $[X] \mapsto [\Hh(X)]$ determine a
contravariant functor $\Hh_k$ from $\rfMm{k}$ to $\Cc$.
\end{thm}
\begin{proof}
We need only show that for $\Lambda_0, \Lambda_1,\Lambda_2 \in
\Obj(\rfMm{k})$ and $[X_i] \in \Hom_{\rfMm{k}}(\Lambda_i,
\Lambda_{i-1})$, there is an isomorphism of
$C^*$-correspondences
\[
\Hh(X_1) \otimes_{C^*(\Lambda_1)} \Hh(X_2)
\cong \Hh(X_1 *_{\Lambda_1^0} X_2).
\]
For $i = 1,2$, let $\Sigma_i$ be a linking graph for $X_i$, and
let $\Sigma_{12}$ be a linking graph for $X_1 *_{\Lambda^0_1}
X_2$. Let $\Gamma$ be the $1$-graph with two edges $a_1, a_2$
and three distinct vertices $v_0 = r(a_1)$, $v_1 = s(a_1) =
r(a_2)$ and $v_2 = s(a_2)$. Then $\Lambda_{v_i} := \Lambda_i$,
$X_{a_i} := X_i$ and $X_{a_1a_2} := X_1 *_{\Lambda^0_1} X_2$
defines a $\Gamma$-system of $k$-morphs ($\theta_{X_1, X_2}$ is
the identity on $X_1 *_{\Lambda^0_1} X_2$). Let $\Sigma$ be a
$\Gamma$-bundle for this system, and let $f : \Sigma \to
\Gamma$ be the bundle map. For $i = 0,1,2$, let $P_i = \sum_{w
\in f^{-1}(v_i) \cap \Sigma^0} s_w \subset MC^*(\Sigma)$.

By applications of the gauge-invariant uniqueness theorem
\cite[Theorem~4.1]{RSY1}, there are canonical isomorphisms
\begin{align*}
C^*(\Sigma_1) &\cong (P_0 + P_1) C^*(\Sigma) (P_0 + P_1) \\
C^*(\Sigma_2) &\cong (P_1 + P_2) C^*(\Sigma) (P_1 + P_2), \text{ and} \\
C^*(\Sigma_{12}) &\cong (P_0 + P_2) C^*(\Sigma) (P_0 + P_2)
\end{align*}
(to establish the third of these isomorphisms, we must slightly
modify the gauge action on $(P_0 + P_2) C^*(\Sigma) (P_0 +
P_2)$). In particular, it follows that
\begin{align*}
\Hh(X_1) &\cong P_0 C^*(\Sigma) P_1 \\
\Hh(X_2) &\cong P_1 C^*(\Sigma) P_2, \text{ and} \\
\Hh(X_1 *_{\Lambda^0_1} X_2) &\cong P_0 C^*(\Sigma) P_2.
\end{align*}
As in Lemma~\ref{lem:inclusions}, the $P_i$ are all full
projections in $MC^*(\Sigma)$, so
\[
\Hh(X_1 *_{\Lambda^0_1} X_2) \cong P_0 C^*(\Sigma) P_2 = (P_0 C^*(\Sigma) P_1) (P_1 C^*(\Sigma) P_2),
\]
and hence, if we identify $\Hh(X_1) \otimes_{C^*(\Lambda_1)}
\Hh(X_2)$ with $(P_0 C^*(\Sigma) P_1) \otimes_{P_1 C^*(\Sigma)
P_1} (P_1 C^*(\Sigma) P_2)$, multiplication in $C^*(\Sigma)$
induces an isomorphism
\[
\Hh(X_1) \otimes_{C^*(\Lambda_1)} \Hh(X_2) \cong
(P_0 C^*(\Sigma) P_1) (P_1 C^*(\Sigma) P_2).
\]
This completes the proof.
\end{proof}

We now present an alternative construction of the
$C^*$-correspondence $\Hh(X)$ (see \cite{Dea} for a similar
construction). In the following, given a $k$-morph $X$, we
denote the point-mass function at $x \in X$ by $\delta_x \in
C_c(X)$. We regard $C_c(X)$ as a right pre-Hilbert
$C_0(\Gamma^0)$ module with pointwise operations.

\begin{prop}\label{prop:correspondence construction}
Let $\Lambda$ and $\Gamma$ be $k$-graphs, and $X$ a \LGmph\
satisfying~\malteseref{}. Let $(\Sigma,i)$ be a linking graph
for $X$, and identify $X$ with $\Sigma^{e_{k+1}}$. Let $\Hh(X)$
be the $C^*$-correspondence obtained from
Proposition~\ref{prp:H(X)}. Then there is an isomorphism of
right-Hilbert $C^*(\Gamma)$-modules
\[
\Hh(X) \cong \overline{C_c(X) \otimes_{C_0(\Gamma^0)} C^*(\Gamma)}
\]
determined by
\begin{equation}\label{eq:module isomorphism}
s_{x i(\alpha)}s^*_{i(\beta)}
\mapsto \delta_x \otimes s_\alpha s^*_{\beta}\qquad\text{ for $x \in \Sigma^{e_{k+1}}$, $\alpha,\beta \in \Gamma$}.
\end{equation}
This isomorphism carries the left action of $s_\lambda \in
C^*(\Lambda)$ on $\Hh(X)$ to
\[
s_\lambda \cdot (\delta_x \otimes s_\alpha s^*_{\beta})
= \delta_{x'} \otimes s_\gamma s_\alpha s^*_\beta \qquad\text{ for $x \in \Sigma^{e_{k+1}}$ and $\alpha,\beta \in \Gamma$}
\]
where $\phi_X(x',\gamma) = (\lambda, x)$.
\end{prop}
\begin{proof}
For the first statement, we just need to check that the
formula~\eqref{eq:module isomorphism} extends to an
inner-product preserving map. For all $x, x' \in
\Sigma^{e_{k+1}}$, $\alpha \in s(x)\Gamma$, $\alpha' \in s(x')
\Gamma$, $\beta \in \Gamma s(\alpha)$ and $\beta' \in \Gamma
s(\alpha')$, we have
\begin{equation}\label{eq:checking module iso}
\langle \delta_x \otimes s_\alpha s^*_{\beta}, \delta_{x'} \otimes s_{\alpha'} s^*_{\beta'} \rangle
 = \langle s_\alpha s^*_\beta, \langle \delta_x, \delta_{x'} \rangle \cdot s_{\alpha'} s^*_{\beta'} \rangle
 =
 \begin{cases}
    s_\beta s^*_\alpha s_{\alpha'} s^*_{\beta'} &\text{ if $x = x'$}\\
    0 &\text{ otherwise.}
 \end{cases}
\end{equation}
Since $s_{x i(\alpha)}^* s_{x'i(\alpha')} = s_{i(\alpha)}^*
(s^*_x s_{x'}) s_{i(\alpha')}$, the third Cuntz-Krieger
relation in $C^*(\Sigma)$ forces
\[
s_{x i(\alpha)}^* s_{x'i(\alpha')} =
\begin{cases}
    s^*_{i(\alpha)} s_{i(\alpha')} &\text{ if $x = x'$}\\
    0 &\text{ otherwise;}
\end{cases}
\]
then equations \eqref{eq:checking module
iso}~and~\eqref{eq:module isomorphism} and the definition of
the inner product on $\Hh(X)$ imply that
\[
\langle \delta_x \otimes s_\alpha s^*_{\beta}, \delta_{x'} \otimes s_{\alpha'} s^*_{\beta'} \rangle
 = \langle s_{x i(\alpha)}s^*_{\beta}, s_{x' i(\alpha')}s^*_{\beta'} \rangle.
\]
By linearity and continuity, it follows that~\eqref{eq:module
isomorphism} is inner-product preserving.

The last assertion follows from a direct calculation.
\end{proof}

We now consider the case where $X$ is a $\Lambda$ endomorph. We
show that the $(k+1)$-graph $C^*$-algebra $C^*(\Lambda \times_X
\NN)$ coincides with the Cuntz-Pimsner algebra $\Oo_{\Hh(X)}$.

Let $X$ be a $\Lambda$ endomorph satisfying~\malteseref{}. Let
$\Lambda \times_X \NN$ be the $T_1$-bundle for the system
induced by $X$ as in
Examples~\ref{ex:skew-products}(\ref{it:endoskewprod}). Because
$T_1$ has just one object, we may simplify the notation of
Remark~\ref{rmk:skewgraphstructure} as follows: there are
injective maps $h_\Lambda : \Lambda \to \Lambda \times_X \NN$
and $h_n : X^{*n} \to \Lambda \times_X \NN$ (where $h_0$ is the
identity map on vertices), and every element of $\Lambda
\times_X \NN$ is of the form $h_n(x)h_\Lambda(\lambda)$ for
some $n \in \NN$, $x \in X^{*n}$ and $\lambda \in \Lambda$.

\begin{thm}\label{thm:CP algebra}
Let $\Lambda$ be a $k$-graph and $X$ a $\Lambda$ endomorph
satisfying~\malteseref{}. Let $\Hh(X)$ be the associated
$C^*(\Lambda)$--\,$C^*(\Lambda)$ correspondence and let
$\Lambda \times_X \NN$ be the $T_1$-bundle associated to $X$
regarded as a $\Lambda$ endomorph. There are a homomorphism
$\pi : C^*(\Lambda) \to C^*(\Lambda \times_X \NN)$ and a linear
map $t : \Hh(X) \to C^*(\Lambda \times_X \NN)$ determined by
\[
\pi(s_\lambda) = s_{h_\Lambda(\lambda)} \qquad\text{and}\qquad
t(s_{x i(\alpha)} s^*_{i(\beta)}) = s_{h_1(x)} s_{h_\Lambda(\alpha)} s^*_{h_\Lambda(\beta)}.
\]
The pair $(t,\pi)$ is a Cuntz-Pimsner covariant representation
of $\Hh(X)$, and the induced $C^*$-homomorphism $t \times \pi :
\Oo_{\Hh(X)} \to C^*(\Lambda \times_X \NN)$ is an isomorphism.
\end{thm}
\begin{proof}
The universal property of $C^*(\Lambda)$ shows that there is a
homomorphism $\pi : C^*(\Lambda) \to C^*(\Lambda \times_X \NN)$
satisfying $\pi(s_\lambda) = s_{h_\Lambda(\lambda)}$ for all
$\lambda \in \Lambda$. This homomorphism $\pi$ is equivariant
for the gauge action on $C^*(\Lambda)$ and the restriction of
the gauge action on $C^*(\Lambda \times_X \NN)$ to the first
$k$ coordinates of $\TT^{k+1}$. Hence an application of the
gauge-invariant uniqueness theorem for $C^*(\Lambda)$
\cite[Theorem~3.4]{KP} shows that $\pi$ is injective.

To see that the formula given for $t$ determines a well-defined
linear map, we will show that for any finite linear combination
of the form $\sum^n_{j=1} a_j s_{x_j i(\alpha_j)}
s^*_{i(\beta_j)}$ in $\Hh(X)$, we have
\[\textstyle
\Big\|\sum^n_{j=1} a_j s_{x_j i(\alpha_j)} s^*_{i(\beta_j)}\Big\|_{\Hh(X)}
=
\Big\|\sum^n_{j=1} a_j
    s_{h_1(x_j)} s_{h_\Lambda(\alpha_j)} s^*_{h_\Lambda(\beta_j)}\Big\|_{C^*(\Lambda \times_X \NN)}.
\]
We have already shown that $\pi$ is injective, so by the
$C^*$-identity for $C^*(\Lambda \times_X \NN)$ and the
definition of the norm on $\Hh(X)$ it suffices to show that for
spanning elements $s_{x i(\alpha)} s^*_{i(\beta)}$ and $s_{x'
i(\alpha')} s^*_{i(\beta')}$ of $\Hh(X)$, we have
\[
\pi(\langle s_{x i(\alpha)} s^*_{i(\beta)}, s_{x' i(\alpha')} s^*_{i(\beta')} \rangle_{\Hh(X)})
= t(s_{x i(\alpha)} s^*_{i(\beta)})^* t(s_{x' i(\alpha')} s^*_{i(\beta')}).
\]
This follows from a routine calculation like~\eqref{eq:checking
module iso} above.

It follows by linearity from the preceding paragraph that
$\pi(\langle \xi, \eta \rangle_{C^*(\Lambda)}) = t(\xi)^*
t(\eta)$ for all $\xi,\eta \in \Hh(X)$. To see that $(t,\pi)$
is a representation, it therefore suffices to show that for a
generator $s_\lambda$ of $C^*(\Lambda)$ and spanning elements
$s_{x i(\alpha)} s^*_{i(\beta)}$ and $s_{x' i(\alpha')}
s^*_{i(\beta')}$ of $\Hh(X)$,
\begin{align*}
\pi(s_\lambda) t(s_{x i(\alpha)} s^*_{i(\beta)})
&= t(s_\lambda \cdot (s_{x i(\alpha)} s^*_{i(\beta)})),\text{ and}\\
t(s_{x i(\alpha)} s^*_{i(\beta)}) \pi(s_\lambda)
&= t((s_{x i(\alpha)} s^*_{i(\beta)}) \cdot s_\lambda).
\end{align*}
One verifies these identities with short calculations using the
definitions of $t$ and $\pi$ and the structure of $\Hh(X)$. We
give the first of these calculations as it is the least
elementary. Fix $s_\lambda$ and $s_{x i(\alpha)}
s^*_{i(\beta)}$ as above, and let $x'$ and $\lambda'$ be the
elements such that $\phi_X(x',\lambda') = (\lambda,x)$. Then
$s_\lambda \cdot (s_{(x i(\alpha))} s^*_{i(\beta)}) = s_{(x'
i(\lambda'\alpha))} s^*_{i(\beta)}$, and we have
\[
\pi(s_\lambda) t(s_{x i(\alpha)} s^*_{i(\beta)})
 = s_{h_\Lambda(\lambda)} s_{h_1(x)} s_{h_\Lambda(\alpha)} s^*_{h_\Lambda(\beta)}
 = s_{h_1(x')} s_{h_\Lambda(\lambda')} s_{h_\Lambda(\alpha)} s^*_{h_\Lambda(\beta)}
 = t(s_\lambda \cdot (s_{(x i(\alpha))} s^*_{i(\beta)})).
\]
To check that $(t,\pi)$ is Cuntz-Pimsner covariant, fix
$\lambda \in \Lambda$, and note that the left action of
$s_\lambda$ on $\Hh(X)$ is implemented by
\[
\sum_{\substack{r(x) = s(\lambda) \\ \phi_X(x',\lambda') = (\lambda, x)}}
\Theta_{s_{x' i(\lambda')}, s_{x}} \in \Kk(\Hh(X)).
\]
Hence if $\varphi$ denotes the homomorphism which implements
the left action, we have
\[
t^{(1)}(\varphi(s_\lambda))
 = \sum_{\substack{r(x) = s(\lambda) \\ \phi_X(x',\lambda') = (\lambda, x)}} t(s_{x' i(\lambda')}) t(s_x)^*
 = \sum_{r(x) = s(\lambda)} \pi(s_{\lambda}) t(s_x) t(s_x)^*,
\]
and this is equal to $\pi(s_\lambda)$ by the fourth
Cuntz-Krieger relation in $C^*(\Lambda \times_X \NN)$.

The restriction of the gauge-action on $C^*(\Lambda \times_X
\NN)$ to the last coordinate in $\TT^{k+1}$ is compatible with
the gauge action on $\Oo_{\Hh(X)}$. The gauge-invariant
uniqueness theorem \cite[Theorem~4.1]{FMR} for $\Oo_{\Hh(X)}$
(see also \cite[Theorem~6.4]{Kat}) now implies that $t \times
\pi$ is injective. It remains only to observe that for
$h_n(x)h_\Lambda(\lambda) \in \Lambda \times_X \NN$, we can
write $x = (x_1, \dots, x_n)$ where each $x_i \in X$, and
$r(x_{i+1}) = s(x_i)$, and then
\[
s_{h_n(x)h_\Lambda(\lambda)} = t(s_{x_1}) \cdots t(s_{x_n}) \pi(s_\lambda)
\]
(if $n = 0$, then $x = r(\lambda) \in I(\Lambda)$, so
$s_{h_n(x)h_\Lambda(\lambda)} = \pi(s_\lambda)$). Hence $t
\times \pi$ is surjective.
\end{proof}

\begin{rmk}
Let $\Sigma$ be a row-finite $(k+1)$-graph with no sources such
that $\Sigma^n v \not= \emptyset$ for all $n \in \NN^{k+1}$ and
$v \in \Sigma^0$. Let $\Sigma^\iota$ and $X$ be the $k$-graph
and $\Sigma^\iota$ endomorph discussed in
Example~\ref{ex:first}(\ref{it:Sigmaiota}). Then $X$
satisfies~\malteseref{}. As in
Example~\ref{ex:skew-products}(\ref{it:alphacrossprod}), we
have that $\Sigma$ is isomorphic to the endomorph skew graph
$\Sigma^\iota \times_X \NN$. Hence Theorem~\ref{thm:CP algebra}
implies that $C^*(\Sigma) \cong \Oo_{\Hh(X)}$; that is the
$C^*$-algebra of $\Sigma$ can be realised as the Cuntz-Pimsner
algebra of a $C^*$-correspondence over $C^*(\Sigma^\iota)$.
\end{rmk}

\begin{rmk}
Let $\Lambda \in \Obj(\rfMm{k})$, let $\alpha$ be an
automorphism of $\Lambda$, and let $X = X(\alpha)$ be the
associated endomorph; clearly $X$ satisfies~\malteseref{}. Let
$\widetilde{\alpha}$ denote the induced automorphism of
$C^*(\Lambda)$. As in
Example~\ref{ex:skew-products}(\ref{it:alphacrossprod}), the
endomorph skew graph $\Lambda \times_X \NN$ is isomorphic to
the crossed product $(k+1)$-graph $\Lambda \times_\alpha \ZZ$
constructed in \cite{FPS}.

In this situation, the bimodule $\Hh(X)$ constructed above is
isomorphic to the bimodule constructed by Pimsner in
\cite[Example~3, p.193]{Pim} with $A = C^*(\Lambda)$ and $\pi =
\widetilde{\alpha}$. We therefore recover the isomorphism
$C^*(\Lambda \times_\alpha \ZZ) \cong C^*(\Lambda)
\times_{\widetilde{\alpha}} \ZZ$ of \cite[Theorem~3.5]{FPS}
from Theorem~\ref{thm:CP algebra} and Pimsner's result.
\end{rmk}

\begin{cor}\label{cor:ste}
Let $\Lambda, \Gamma$ be $k$-graphs, let $R$ be a \LGmph\ and
$S$ a \morph{\Gamma}{\Lambda}. Suppose $R$ and $S$
satisfy~\malteseref{}. Let $X$ be the $\Lambda$ endomorph $S
*_{\Gamma^0} R$, and let $Y$ be the $\Gamma$-endomorph $R
*_{\Lambda^0} S$. Then the $(k+1)$-graph algebras $C^*(\Lambda
\times_X \NN)$ and $C^*(\Gamma \times_Y \NN)$ are Morita
equivalent.

In particular the above Morita equivalence holds if $X =
\sideset{_p}{_q}{\opX}$ and $Y = X_q *_{\Lambda^0}
\sideset{_p}{}{\opX}$ where $p,q : \Gamma \to \Lambda$ are
finite coverings of row-finite $k$-graphs.
\end{cor}
\begin{proof}
Theorem~\ref{thm:functoriality} implies that $\Hh(X) \cong
\Hh(S) \otimes_{C^*(\Gamma)} \Hh(R)$, and $\Hh(Y) \cong \Hh(R)
\otimes_{C^*(\Lambda)} \Hh(S)$. By
Theorem~\ref{prop:correspondence construction}, $E := \Hh(X)$,
and $F := \Hh(Y)$ satisfy the hypotheses of
\cite[Theorem~3.14]{MPT} (the $C^*$-correspondences $\Hh(R)$
and $\Hh(S)$ implement the elementary strong shift
equivalence). Hence $\Oo_{\Hh(X)}$ and $\Oo_{\Hh(Y)}$ are
Morita equivalent, and our result follows from
Theorem~\ref{thm:CP algebra}.
\end{proof}

\begin{rmk}
When $k = 0$, the first statement of the above Corollary
reduces to Bates' results on shift-equivalence for $1$-graphs
\cite{Bates}.
\end{rmk}

\begin{example}
For $n \in \NN\setminus\{0\}$, let $D_n$ be the directed graph
with $n$ vertices $\{v_0, \dots, v_{n-1}\}$ and edges $\{x_i,
y_i : 0 \le i \le n-1\}$ where $r(x_i) = v_i = s(y_i)$ and
$s(x_i) = v_{i+1} = r(y_i)$ (see \cite[Section~6.2]{KPS}). In
particular, $D_1$ is equal to the bouquet of two loops whose
$C^*$-algebra is canonically isomorphic to $\Oo_2$. We will
consider a $D_1$ endomorph $\sideset{_p}{_q}{\opX}$ constructed
from coverings $p,q : D_2 \to D_1$. To reduce confusion, we
will denote $x_0, y_0 \in D_1$ by $x$ and $y$.

The covering maps $p,q : D_2 \to D_1$ are defined as follows
\begin{align*}
&&&&&&p(x_i) &= x, & p(y_i) &= y, &&&&&&\\
&&&&&&q(x_0) &= x, & q(y_0) &= y, &&&&&&\\
&&&&&&q(x_1) &= y, & q(y_1) &= x. &&&&&&
\end{align*}
Construct $X = \sideset{_p}{_q}{\opX}$ as in
Example~\ref{ex:first}(\ref{it:pXq}) (so, as a set, $X =
D_2^0$). The endomorph skew-graph $\Lambda = D_1 \times_X \NN$
is a $2$-graph whose skeleton and factorisation rules can be
described as follows: $\Lambda^0 = v$, $\Lambda^{e_1} =
\{x,y\}$, $\Lambda^{e_2} = X = \{v_0, v_1\}$, and
\[
v_0 x = y v_1, \quad v_0 y = x v_1, \quad v_1 x = x v_0, \quad\text{and}\quad v_1 y = y v_0.
\]
Results of \cite{APS, RobSt, RobSt2} can be used to see that
$C^*(\Lambda)$ is a Kirchberg algebra (the details appear in an
unpublished manuscript of D. Robertson \cite{DRNotes}) and has
trivial $K$-theory. Hence $C^*(\Lambda)$ is isomorphic $\Oo_2$
by the Kirchberg-Phillips theorem.

Using Proposition~\ref{prop:correspondence construction}, we
see that $\Hh(X) \cong \CC^2 \otimes_{\CC} \Oo_2$ as a
right-Hilbert $\Oo_2$-module (where $\CC^2$ is the
two-dimensional Hilbert space with orthonormal basis
$\{\delta_0, \delta_1\}$). The left action of $\Oo_2 =
C^*(\{S_0, S_1\})$ is determined by
\begin{align*}
&&&& S_0 \cdot (\delta_0 \otimes_\CC 1_{\Oo_2}) &= \delta_1 \otimes_{\CC} S_0, &
S_0 \cdot (\delta_1 \otimes_\CC 1_{\Oo_2}) &= \delta_0 \otimes_{\CC} S_1, &&&& \\
&&&& S_1 \cdot (\delta_0 \otimes_\CC 1_{\Oo_2}) &= \delta_1 \otimes_{\CC} S_1, &
S_1 \cdot (\delta_1 \otimes_\CC 1_{\Oo_2}) &= \delta_0 \otimes_{\CC} S_0. &&&&
\end{align*}
By Theorem~\ref{thm:CP algebra}, the Cuntz-Pimsner algebra
$\Oo_{\Hh(X)}$ of this $C^*$-correspondence is isomorphic to
$C^*(\Lambda)$ which, as we saw above, is isomorphic to
$\Oo_2$.

Let $Y = X_q *_{\Lambda^0} \sideset{_p}{}{\opX}$. Note that $Y$
satisfies~\malteseref{}. The endomorph skew-graph $D_2 \times_Y
\NN$ has skeleton
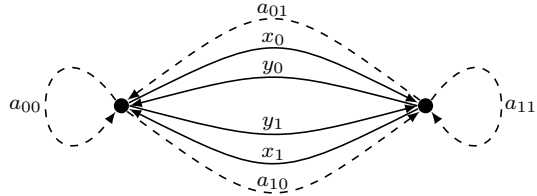
\[
\begin{tikzpicture}[scale=2]
    \node[inner sep=2pt, fill=black, circle] (v) at (-1,0) {};
    \node[inner sep=2pt, fill=black, circle] (w) at (1,0) {};
    \draw[style=semithick, -latex] (w.west) .. controls (0, 0.25) and (0,0.25) .. (v.east)
        node[pos=0.5, anchor=south, inner sep=1pt]{$\scriptstyle y_0$};
    \draw[style=semithick, -latex] (v.22) .. controls (0, 0.5) and (0,0.5) .. (w.158)
        node[pos=0.5, anchor=south, inner sep=1pt]{$\scriptstyle x_0$};
    \draw[style=semithick, -latex] (v.east) .. controls (0, -0.25) and (0,-0.25) .. (w.west)
        node[pos=0.5, anchor=south, inner sep=1pt]{$\scriptstyle y_1$};
    \draw[style=semithick, -latex] (w.202) .. controls (0, -0.5) and (0,-0.5) .. (v.338)
        node[pos=0.5, anchor=south, inner sep=1pt]{$\scriptstyle x_1$};
    \draw[style=dashed, style=semithick, -latex] (w.north west) .. controls (0, 0.75) and (0,0.75) .. (v.north east)
        node[pos=0.5, anchor=south, inner sep=1pt]{$\scriptstyle a_{01}$};
    \draw[style=dashed, style=semithick, -latex] (v.south east) .. controls (0, -0.75) and (0,-0.75) .. (w.south west)
        node[pos=0.5, anchor=south, inner sep=1pt]{$\scriptstyle a_{10}$};
    \draw[style=dashed, style=semithick, -latex] (w.north east) .. controls (1.3, 0.4) and (1.5, 0.25) .. (1.5,0)
                                                          .. controls (1.5, -0.25) and (1.3, -0.4) .. (w.south east);
    \node[anchor=east, inner sep=1pt] at (-1.5,0) {$\scriptstyle a_{00}$};
    \draw[style=dashed, style=semithick, -latex] (v.north west) .. controls (-1.3, 0.4) and (-1.5, 0.25) .. (-1.5,0)
                                                          .. controls (-1.5, -0.25) and (-1.3, -0.4) .. (v.south west);
    \node[anchor=west, inner sep=1pt] at (1.5,0) {$\scriptstyle a_{11}$};
\end{tikzpicture}
\]
with factorisation rules
\begin{align*}
a_{00}y_0 &= y_0 a_{11} & a_{00} x_1 &= x_1 a_{11} & a_{10} y_0 &= x_0 a_{01} & a_{10} x_1 = y_1 a_{01} \\
a_{11}x_0 &= y_1 a_{01} & a_{11} y_1 &= x_0 a_{01} & a_{01} x_0 &= x_0 a_{10} & a_{01} y_1 = y_1 a_{10}
\end{align*}
Corollary~\ref{cor:ste} shows that $C^*(D_2 \times_Y \NN)$ is
Morita equivalent to $C^*(\Lambda) \cong \Oo_2$, and as it is
also unital, it is in fact isomorphic to $\Oo_2$.
\end{example}

\end{document}